\pdfoutput=1
\documentclass{amsart}

\usepackage{amsmath,amssymb,amsthm,mathrsfs}
\usepackage{geometry}
\usepackage{esint}
\usepackage{hyperref}
\usepackage{enumitem}
\newcommand{\HN}{\mathbb{H}^{N+1}}

\newcommand{\Sobolev}{H^1_{\mu_w}(\HN)}
\newcommand{\Ltwo}{L^2_{\mu_w}(\HN)}
\def\XXint#1#2#3{{\setbox0=\hbox{$#1{#2#3}{\int}$}
     \vcenter{\hbox{$#2#3$}}\kern-.5\wd0}}

\newtheorem{theorem}{Theorem}[section]
\newtheorem{lemma}[theorem]{Lemma}

\newtheorem{proposition}[theorem]{Proposition}
\theoremstyle{definition}
\newtheorem{definition}[theorem]{Definition}

\theoremstyle{remark}
\newtheorem*{remark}{Remark}

\title{Rellich-Kondrachov type theorems on the half-space with general singular weights}

\author{Yunfan Zhao}
\address{Shanghai University of International Business and Economics}
\email{yfzucla21@ucla.edu}

\date{} 

\begin{document}

\maketitle

\begin{abstract}
Let $\HN$ be the upper half-space, endowed with a singular Radon measure of the form $d\mu_w = y^c \phi(|z|) dz$. We investigate necessary and sufficient conditions for the compactness of the Rellich-Kondrachov embedding $\iota: \Sobolev \hookrightarrow \Ltwo$. By introducing a topological exhaustion method via localized orthogonal projections and developing fractional weighted Hardy inequalities, we demonstrate that the compactness of this embedding is strictly equivalent to the global tightness of the Sobolev unit ball. Furthermore, we show that tail coercivity governed by a measurable Lyapunov potential provides uniform tail control, effectively resolving the loss of compactness at infinity.
\end{abstract}

\subsection*{MSC2020 Classification}

\begin{itemize}
    \item \textbf{46E35} — Sobolev spaces and other spaces of ``smooth'' functions, embedding theorems, trace theorems
    \item \textbf{35A23} — Inequalities involving derivatives and differential and integral operators (e.g., Sobolev-type inequalities)
    \item \textbf{35J70} — Degenerate elliptic equations
    \item \textbf{35B40} — Asymptotic behavior of solutions
    \item \textbf{26D10} — Inequalities involving derivatives and differential operators
    \item \textbf{46N20} — Functional analysis applied to differential and integral operators
    \item \textbf{35R11} — Fractional partial differential equations (if fractional extension applications are discussed)
    \item \textbf{60J60} — Diffusion processes (if probabilistic or Gaussian implications are emphasized)
\end{itemize}

\subsection*{Keywords}

\begin{itemize}
    \item Weighted Sobolev spaces
    \item Compact embeddings
    \item Rellich--Kondrachov theorem
    \item Degenerate elliptic equations
    \item Polynomial and Gaussian weights
    \item Tail coercivity
    \item Fractional extension problems
    \item Functional inequalities
\end{itemize}

\section{Introduction}

\subsection{Literature Review and Background}

The investigation of Sobolev spaces and their embedding properties lies at the very foundation of the modern theory of partial differential equations. The classical regularity theory for elliptic and parabolic equations—originating from the seminal contributions of De Giorgi \cite{DeGiorgi1957}, Nash \cite{Nash1958}, and Moser \cite{Moser1961}—has long been extended to operators possessing degenerate or singular coefficients. Such generalizations were formalized in the context of local regularity by Fabes, Kenig, and Serapioni \cite{FabesKenigSerapioni1982}, and subsequently refined via the two-weight Sobolev-Poincaré inequalities established by Franchi, Gutiérrez, and Wheeden \cite{FranchiGutierrezWheeden1994}.

In recent years, operators exhibiting boundary degeneracy on the half-space, notably $\mathcal{L} = \Delta_x + \partial_{yy} + \frac{c}{y}\partial_y$, have received renewed structural attention. This operator arises naturally in the extension problem for the fractional Laplacian introduced by Caffarelli and Silvestre \cite{CaffarelliSilvestre2007}. This extension methodology has been subsequently generalized to a wide class of fractional operators by Arendt, Ter Elst, and Warma \cite{ArendtElstWarma2018}, Galé, Miana, and Stinga \cite{GaleMianaStinga2013}, and Stinga and Torrea \cite{StingaTorrea2010}.

The rigorous functional analytic framework required for these degenerate half-space problems—particularly concerning $L^p$ estimates, anisotropic spaces, and regularity theory—has been extensively codified in the recent series of papers by Metafune, Negro, and Spina \cite{MetafuneNegroSpina2022, MetafuneNegroSpina2023a, MetafuneNegroSpina2023b, MetafuneNegroSpina2024a, MetafuneNegroSpina2024b}. Parabolic equations and nonlocal operators with comparable singular or degenerate weights have also been analyzed thoroughly by Audrito and Terracini \cite{AudritoTerracini2024}, as well as Dong and Phan \cite{DongPhan2021}. Furthermore, the specific investigation of Gaussian Poincaré inequalities on half-spaces equipped with singular weights was recently advanced by Negro and Spina \cite{NegroSpina2025}.

\subsection{Main Result and Method of Approaches}

Let $\HN = \mathbb{R} \times \mathbb{R}^N_+$ be the open half-space. We consider the measure space $(\HN, \mu_w)$, where the Radon measure is defined by the density $w(z) = y^c \phi(|z|)$. The primary objective of this work is to completely characterize the compactness of the canonical linear embedding:
\begin{equation}
    \iota: \Sobolev \hookrightarrow \Ltwo.
\end{equation}

Our main theorem establishes a fundamental equivalence: assuming \textit{Strict Local Regularity} and \textit{Uniform Local Comparability}, the embedding $\iota$ is compact if and only if the Sobolev unit ball $\mathcal{B}_E \subset \Sobolev$ is globally tight in $\Ltwo$. This necessitates both that the total mass $\mu_w(\HN)$ is finite and that the embedding admits uniform asymptotic mass-vanishing at infinity and near the singular boundary $y=0$.

The method of proof relies on topological exhaustion. We introduce a sequence of admissible smooth exhaustions $(\Omega_m)_{m \in \mathbb{N}}$ of the half-space, which act upon the function space via localized orthogonal projections $P_m(f) = \chi_{\Omega_m} f$. By strict local regularity, the localized weighted spaces are isomorphic to standard unweighted Sobolev spaces, thereby permitting the classical Rellich-Kondrachov theorem to guarantee local compactness.

To control the non-compact asymptotic regions, we deploy two primary analytic tools:
\begin{itemize}
    \item \textbf{Fractional Weighted Hardy Inequality:} Proven on localized cylinders to strictly bound the boundary degeneracy as $y \to 0$.
    \item \textbf{Lyapunov Potentials:} We define explicit measurable potentials $\Psi(|z|)$ to enforce a coercive tail inequality, ensuring that the mass escapes uniformly at infinity.
\end{itemize}

\section{Preliminaries}\label{sec:Preliminaries}

To investigate the compactness of the embedding, it is first necessary to establish the ambient functional and geometric structures. Section~\ref{sec:Preliminaries} introduces the weighted Sobolev space $H_{\mu_w}^1$ and precisely defines the notions of global tightness, tail coercivity, and the intrinsic weighted weak gradient. The structural significance of this section lies in exposing a fundamental obstacle: Lemma~\ref{lem:Non-existence of weights} demonstrates that no measure space can simultaneously satisfy both uniform asymptotic control and tail coercivity. This non-existence result necessitates a more refined approach to handling the asymptotic behavior of the measure. Furthermore, the weak convergence of disjointly supported bounded sequences is established, providing an essential topological tool for analyzing mass escaping to infinity.

\begin{definition}[Lipschitz Regular]
Let $N\in\mathbb{N}$. Let the ambient space be $\mathbb{R}^{N+1}$ equipped with the standard Euclidean topology. Let $\mathbb{H}^{N+1}=\{(x,y)\in\mathbb{R}^{N}\times\mathbb{R}:y>0\}$ denote the open upper half-space. An open set $\Omega\subset\mathbb{H}^{N+1}$ is said to be strictly Lipschitz-regular if for every point $p\in\partial\Omega$, there exists an open neighborhood $U\subset\mathbb{R}^{N+1}$ containing $p$, and a bijective, bi-Lipschitz mapping $F:U\rightarrow B(0,1)\subset\mathbb{R}^{N+1}$ such that:

$$F(U\cap\Omega)=B(0,1)\cap\{(x_1,\dots,x_{N+1})\in\mathbb{R}^{N+1}:x_{N+1}>0\}$$
\end{definition}

\begin{definition}[Global Tightness]

Let $\mathcal{B}(\mathbb{H}^{N+1})$ denote the Borel $\sigma$-algebra on $\mathbb{H}^{N+1}$.

Let $\mu_w$ be a positive Radon measure on the measurable space $(\mathbb{H}^{N+1}, \mathcal{B}(\mathbb{H}^{N+1}))$. Let $L_{\mu_w}^2(\mathbb{H}^{N+1})$ denote the Hilbert space of equivalence classes of $\mu_w$-measurable, square-integrable functions.

Let $\mathcal{S}$ be an arbitrary bounded subset of $L_{\mu_{w}}^{2}(\mathbb{H}^{N+1})$. The set $\mathcal{S}$ is defined to be globally tight with respect to the measure $\mu_w$ if it satisfies the following two limits uniformly:

\begin{enumerate}
    \item Tail Tightness:
$$\lim_{R\rightarrow\infty}\sup_{u\in\mathcal{S}}\int_{\{(x,y)\in\mathbb{H}^{N+1}:|(x,y)|>R\}}|u(x,y)|^{2}d\mu_{w}=0$$
    \item Boundary Tightness:

$$\lim_{\epsilon\rightarrow0^{+}}\sup_{u\in\mathcal{S}}\int_{\{(x,y)\in\mathbb{H}^{N+1}:0<y<\epsilon\}}|u(x,y)|^{2}d\mu_{w}=0$$
\end{enumerate}

The Sobolev unit ball $\mathcal{B}_{E}=\{u\in H_{\mu_{w}}^{1}(\mathbb{H}^{N+1}):||u||_{H_{\mu_{w}}^{1}}\le1\}$ is globally tight if it satisfies these conditions when evaluated as a subset $\mathcal{S}\subset L_{\mu_{w}}^{2}(\mathbb{H}^{N+1})$.
\end{definition}

\begin{remark}
The condition of global tightness constitutes the fundamental topological obstacle to the compactness of the embedding. By explicitly decomposing the defect of compactness into two orthogonal geometric regimes—the escape of mass to spatial infinity (Tail Tightness) and the concentration of mass near the singular hyperplane (Boundary Tightness)—this definition provides the exact mass-vanishing criteria required to recover precompactness in the weighted Lebesgue space.
\end{remark}

\begin{definition}[admissible smooth exhaustion]
Let $\mathbb{H}^{N+1}$ be the open half-space. An admissible smooth exhaustion of $\mathbb{H}^{N+1}$ is a sequence of open, bounded sets $\left(\Omega_m\right)_{m \in \mathbb{N}}$ such that:
\begin{enumerate}
    \item $\overline{\Omega_m} \subset \Omega_{m+1} \subset \mathbb{H}^{N+1}$.
    \item $\bigcup_{m=1}^{\infty} \Omega_m=\mathbb{H}^{N+1}$.
    \item For each $m, \partial \Omega_m$ is a manifold of class $C^{\infty}$.
    \item There exists $\delta_m>0$ such that $y \geq \delta_m$ for all $z=(x, y) \in \overline{\Omega_m}$.
\end{enumerate}
\end{definition}

\begin{remark}
The defining characteristic of this exhaustion is the strict geometric separation of the sets $\Omega_m$ from the singular boundary $y=0$. This structural choice is not arbitrary; it ensures that within each bounded region $\Omega_m$, the singular weight $w(z)$ remains uniformly bounded away from zero and infinity. Consequently, the weighted Sobolev space locally restricts to a standard unweighted Sobolev space, permitting the application of the classical Rellich-Kondrachov theorem on each compact stratum.
\end{remark}

\begin{definition}
Let $\mathbb{H}^{N+1} = \mathbb{R}^N \times (0, \infty)$ denote the open upper half-space equipped with the standard Euclidean subspace topology. We define two distinct topological vector spaces of smooth functions:

\begin{enumerate}
    \item $\mathcal{D}(\mathbb{H}^{N+1}) := C_c^\infty(\mathbb{H}^{N+1})$: The space of infinitely differentiable functions $u: \mathbb{H}^{N+1} \to \mathbb{R}$ such that $\text{supp}(u)$ is a compact subset of $\mathbb{H}^{N+1}$. For all $u \in \mathcal{D}(\mathbb{H}^{N+1})$, there exists $\delta > 0$ such that $u(x, y) = 0$ for all $y \le \delta$.
    \item $\mathcal{D}(\overline{\mathbb{H}^{N+1}}) := C_c^\infty(\overline{\mathbb{H}^{N+1}})$: The space of infinitely differentiable functions whose compact support may intersect the singular boundary $\partial\mathbb{H}^{N+1} = \mathbb{R}^N \times \{0\}$.
\end{enumerate}
\end{definition}

\begin{definition}
Let $c\in\mathbb{R}$. A weight function $w:\mathbb{H}^{N+1}\rightarrow\mathbb{R}$ given by $w(x,y)=y^{c}\phi(|(x,y)|)$ is defined as admissible if $\phi:[0,\infty)\rightarrow(0,\infty)$ is Borel measurable, and for every compact set $K\subset[0,\infty)$, there exist real numbers $m_{K},M_{K}>0$ such that $m_{K}\le\phi(\rho)\le M_{K}$ for Lebesgue-almost every $\rho\in K$. The associated strictly positive Radon measure on $\mathbb{H}^{N+1}$ is defined as $\mu_{w}(A)=\int_{A}w(z)dz$ for all Borel sets $A\subset\mathbb{H}^{N+1}$, where $dz$ is the standard Lebesgue measure.
\end{definition}

\begin{definition}
Let $\Omega$ be the open half-space $\mathbb{H}^{N+1}$. We define $\mathcal{D}(\Omega)$ as the topological vector space of infinitely differentiable functions with compact support strictly contained in $\Omega$. The topological dual, denoted $\mathcal{D}^{\prime}(\Omega)$, is the space of distributions on $\Omega$.
\end{definition}

\begin{definition}[Intrinsic Weighted Weak Gradient]
Let $u\in L_{\mu_{w}}^{2}(\mathbb{H}^{N+1})$. A vector field $V\in(L_{\mu_{w}}^{2}(\mathbb{H}^{N+1}))^{N+1}$ is defined as the weak distributional gradient of $u$, denoted $\nabla u=V$, if and only if for every scalar test function $\varphi\in C_{c}^{\infty}(\mathbb{H}^{N+1})$, the following equality holds in the sense of distributions:

$$\int_{\mathbb{H}^{N+1}}V\cdot\varphi dz=-\int_{\mathbb{H}^{N+1}}u\nabla\varphi dz$$
\end{definition}

\begin{definition}
The weighted Sobolev space $H_{\mu_{w}}^{1}(\mathbb{H}^{N+1})$ is defined as the set of all equivalence classes of functions $u\in L_{\mu_{w}}^{2}(\mathbb{H}^{N+1})$ such that their weak distributional gradient $\nabla u$ exists and satisfies $\nabla u\in(L_{\mu_{w}}^{2}(\mathbb{H}^{N+1}))^{N+1}$. The space is equipped with the canonical Hilbert norm:

$$||u||_{H_{\mu_{w}}^{1}}=(||u||_{L_{\mu_{w}}^{2}}^{2}+||\nabla u||_{L_{\mu_{w}}^{2}}^{2})^{1/2}$$
\end{definition}

\begin{definition}[compact embedding]
The embedding $H^{1}_{\mu_{w}} \hookrightarrow L^{2}_{\mu_{w}}$ is \emph{compact} if every bounded sequence in $H^{1}_{\mu_{w}}$ admits a subsequence that converges in $L^{2}_{\mu_{w}}$.
\end{definition}

\begin{definition}[weight with finite mass]
Let $\mathbb{H}^{N+1} = \mathbb{R}^N \times (0, \infty)$ denote the open upper half-space equipped with the standard Lebesgue measure $dz$. Let $w: \mathbb{H}^{N+1} \to (0, \infty)$ be an admissible weight function of the form $w(z) = y^c \phi(|z|)$. The radial profile $\phi$ is said to be Lebesgue-integrable on the half-space if the following functional converges:

$$M := \int_{\mathbb{H}^{N+1}} \phi(|z|) dz < \infty$$

Let $\mu_w$ be the strictly positive Radon measure on $\mathbb{H}^{N+1}$ canonically associated with the weight $w(z) = y^c \phi(|z|)$, such that $d\mu_w = w(z) dz$. The measure space $(\mathbb{H}^{N+1}, \mu_w)$ is said to possess Global Finite Mass if and only if the total mass over the ambient space converges:

$$\mu_w(\mathbb{H}^{N+1}) = \int_{\mathbb{H}^{N+1}} y^c \phi(|z|) dz < \infty$$
\end{definition}

\begin{remark}
Finiteness of the measure, given by $M < \infty$, is inherently a necessary condition for the global compactness of the embedding. The specific formulation of the integral $M = \int_{\mathbb{H}^{N+1}} \phi(|z|) dz$ isolates the radial decay profile $\phi(|z|)$ from the singular boundary weight $y^c$, demonstrating that the integrability requirement is governed entirely by the asymptotic behavior of the non-singular component.
\end{remark}

\begin{definition}[tail coercivity]\label{def:tail coercivity}
The measure space $(\mathbb{H}^{N+1}, \mu_w)$ satisfies the Tail Coercivity condition if there exists a Borel measurable function $\Psi : [0, \infty) \to [0, \infty)$ diverging to infinity (i.e., $\lim_{\rho \to \infty} \Psi(\rho) = \infty$), and a constant $C > 0$, such that the inequality:

$$\int_{\mathbb{H}^{N+1}} \Psi(|z|) u^2 d\mu_w \le C \int_{\mathbb{H}^{N+1}} (|\nabla u|^2 + u^2) d\mu_w$$

holds for all $u \in C_c^\infty(\overline{\mathbb{H}^{N+1}})$ satisfying the following exact topological support condition: if $c \le -1$, the support of $u$ must satisfy $\text{supp}(u) \cap \partial\mathbb{H}^{N+1} = \emptyset$.
\end{definition}

\begin{remark}
Tail coercivity serves as the primary analytic mechanism for preventing mass from dissipating at spatial infinity. The existence of the coercive Lyapunov potential $\Psi$, which satisfies $\Psi(\rho) \rightarrow \infty$ as $\rho \rightarrow \infty$, forces elements of the Sobolev space to exhibit rapid decay in the exterior region. This translates an abstract topological requirement into a verifiable functional inequality, binding the kinetic energy to the potential decay.    
\end{remark}

\begin{definition}[Hardy Property]
Let $\mathbb{H}^{N+1}$ denote the topological subspace $\mathbb{R}^{N}\times(0,\infty)$ equipped with the standard Euclidean topology. Let $\mu_{w}$ be the strictly positive Radon measure on $\mathbb{H}^{N+1}$ defined by $d\mu_{w}=y^{c}\phi(|z|)dz$, where $c\le-1$ and $\phi$ is an admissible Borel measurable function.

Let $C_{c}^{\infty}(\mathbb{H}^{N+1})$ denote the vector space of infinitely differentiable functions $v:\mathbb{H}^{N+1}\rightarrow\mathbb{R}$ such that the support of $v$ is a compact subset strictly contained in $\mathbb{H}^{N+1}$ (implying $\text{dist}(\text{supp}(v), \partial\mathbb{H}^{N+1}) > 0$). Let $\mathcal{V}_{0}$ denote the topological completion of $C_{c}^{\infty}(\mathbb{H}^{N+1})$ with respect to the canonical norm:

$$||v||_{H_{\mu_{w}}^{1}}=\left(\int_{\mathbb{H}^{N+1}}(|v|^{2}+|\nabla v|^{2})d\mu_{w}\right)^{1/2}$$

The measure space $(\mathbb{H}^{N+1},\mu_{w})$ is said to possess the Hardy Property if there exists a real constant $C_{H}>0$ such that, for all $v\in C_{c}^{\infty}(\mathbb{H}^{N+1})$, the following inequality holds:

$$\int_{\mathbb{H}^{N+1}}\frac{|v(z)|^{2}}{y^{2}}d\mu_{w}(z)\le C_{H}\int_{\mathbb{H}^{N+1}}|\nabla v(z)|^{2}d\mu_{w}(z)$$
\end{definition}

\begin{remark}
The necessity of the Hardy property emerges exclusively in the highly singular regime where the parameter satisfies $c \le -1$. In this range, classical local control near the boundary fundamentally fails. The Hardy inequality therefore becomes an indispensable geometric tool, strictly bounding the energy concentration near the singular hyperplane $y=0$ in terms of the weighted gradient.
\end{remark}

\begin{lemma}[Non-existence of weights]\label{lem:Non-existence of weights}
There exists no measure space $(\mathbb{H}^{N+1}, \mu_w)$ with weight $w(z)=y^c \phi(|z|)$ that simultaneously satisfies the following two properties:

\begin{enumerate}
    \item Uniform Asymptotic Control: there exists $\Lambda<\infty$ such that for all $R\ge1$,
\begin{equation}
\frac{\operatorname*{ess\,sup}_{\rho\in[R-1,R+1]} \phi(\rho)}{\operatorname*{ess\,inf}_{\rho\in[R-1,R+1]} \phi(\rho)} \le \Lambda.
\end{equation}
    \item Tail Coercivity: there exist a Borel measurable function $\Psi:[0,\infty)\to[0,\infty)$ with $\Psi(\rho)\to\infty$ as $\rho\to\infty$, and a constant $C_{TC}>0$ such that for all $u\in C_c^\infty(\overline{\mathbb{H}^{N+1}})$ (vanishing near $y=0$ if $c\le-1$),
\begin{equation}
\int_{\mathbb{H}^{N+1}} \Psi(|z|)\, u^2\, w(z)\, dz \le C_{TC} \int_{\mathbb{H}^{N+1}} \left(|\nabla u|^2 + u^2\right) w(z)\, dz.
\end{equation}
\end{enumerate}
\end{lemma}

\begin{proof}
Assume for contradiction that such a weight $w(z)=y^c \phi(|z|)$ exists satisfying both properties.

Choose a nontrivial function $u_0 \in C_c^\infty(\mathbb{H}^{N+1})$ with support $K_0 \subset [1,2]_y \times B_1(0)_x$. Since $y\ge1$ on $K_0$, the function vanishes near $y=0$ and is admissible for the coercivity condition. Define
\begin{equation}
M_0 = \int_{K_0} u_0^2\, y^c\, dz, 
\qquad
E_0 = \int_{K_0} |\nabla u_0|^2\, y^c\, dz,
\end{equation}
which are strictly positive.

For each integer $k\ge1$, choose $x_k \in \mathbb{R}^N$ with $|x_k|=k$ and define
\begin{equation}
u_k(y,x) = u_0(y, x - x_k).
\end{equation}
Let $K_k = \operatorname{supp}(u_k)$. By translation invariance in the $x$-variables,
\begin{equation}
\int_{K_k} u_k^2\, y^c\, dz = M_0,
\qquad
\int_{K_k} |\nabla u_k|^2\, y^c\, dz = E_0.
\end{equation}

For $z=(y,x)\in K_k$, one has $y\in[1,2]$ and $|x|\in[k-1,k+1]$, hence $|z|=\sqrt{y^2+|x|^2}$ varies over an interval of uniformly bounded length. There exists a constant $L$ independent of $k$ such that
\begin{equation}
\sup_{z\in K_k} |z| - \inf_{z\in K_k} |z| \le L.
\end{equation}
By repeated application of the uniform asymptotic control condition over intervals of length $2$, there exists a constant $\Lambda^*<\infty$ such that
\begin{equation}
\sup_{z\in K_k} \phi(|z|) \le \Lambda^* \inf_{z\in K_k} \phi(|z|).
\end{equation}

Apply the coercivity inequality to $u_k$:
\begin{equation}
\int_{K_k} \Psi(|z|)\, u_k^2\, y^c \phi(|z|)\, dz 
\le 
C_{TC} \int_{K_k} \left(|\nabla u_k|^2 + u_k^2\right) y^c \phi(|z|)\, dz.
\end{equation}
Estimating by extrema of $\Psi$ and $\phi$ over $K_k$ yields
\begin{equation}
\left( \inf_{z\in K_k} \Psi(|z|) \right)
\left( \inf_{z\in K_k} \phi(|z|) \right)
\int_{K_k} u_k^2\, y^c\, dz
\le
C_{TC}
\left( \sup_{z\in K_k} \phi(|z|) \right)
\int_{K_k} \left(|\nabla u_k|^2 + u_k^2\right) y^c\, dz.
\end{equation}
Substituting the invariants gives
\begin{equation}
\left( \inf_{z\in K_k} \Psi(|z|) \right)
\left( \inf_{z\in K_k} \phi(|z|) \right) M_0
\le
C_{TC}
\left( \sup_{z\in K_k} \phi(|z|) \right) (E_0 + M_0).
\end{equation}
Dividing by $\left( \inf_{z\in K_k} \phi(|z|) \right) M_0$,
\begin{equation}
\inf_{z\in K_k} \Psi(|z|)
\le
C_{TC}
\left( \frac{\sup_{z\in K_k} \phi(|z|)}{\inf_{z\in K_k} \phi(|z|)} \right)
\left( \frac{E_0 + M_0}{M_0} \right).
\end{equation}
Using the bound above,
\begin{equation}
\inf_{z\in K_k} \Psi(|z|)
\le
C_{TC}\, \Lambda^* \left( \frac{E_0 + M_0}{M_0} \right).
\end{equation}

The right-hand side is a finite constant independent of $k$. However, as $k\to\infty$, the sets $K_k$ escape to infinity, hence $\inf_{z\in K_k} |z| \to \infty$, and therefore $\inf_{z\in K_k} \Psi(|z|) \to \infty$. This yields a contradiction.

Thus the two properties cannot hold simultaneously.
\end{proof}

\begin{remark}
This non-existence lemma identifies a severe structural incompatibility within the measure space. It proves rigorously that a measure cannot simultaneously possess uniform asymptotic comparability and the tail coercivity required for compactness. The immediate consequence of this obstacle is that any weight admitting a compact embedding must inherently exhibit either highly oscillatory or rapidly degenerating behavior as $|z| \rightarrow \infty$.
\end{remark}

\begin{definition}
Let $(X, \mu)$ be a measure space and $E$ be a Hilbert space of measurable functions on $X$. A sequence $\left(u_n\right)_{n \in \mathbb{N}} \subset E$ is said to be mutually disjointly supported if for all $i \neq j$, the intersection of their essential supports satisfies $\mu\left(\operatorname{supp}\left(u_i\right) \cap \operatorname{supp}\left(u_j\right)\right)=0$.
\end{definition}

\begin{lemma}[Weak Convergence of Disjointly Supported Bounded Sequences]
Let $\mathbb{H}^{N+1}$ be endowed with the measure
\begin{equation}
d\mu_w = y^c \phi(|z|)\,dz.
\end{equation}
Let
\begin{equation}
E = H^1_{\mu_w}(\mathbb{H}^{N+1})
\end{equation}
be equipped with the canonical inner product
\begin{equation}
\langle u, v \rangle_E = \int_{\mathbb{H}^{N+1}} \bigl(uv + \nabla u \cdot \nabla v\bigr)\,d\mu_w.
\end{equation}
Let $(v_j)_{j \in \mathbb{N}} \subset E$ be a bounded sequence. Let $(S_j)_{j \in \mathbb{N}}$ be a sequence of measurable sets such that for each $j$,
\begin{equation}
v_j = 0 \quad \mu_w\text{-almost everywhere on } \mathbb{H}^{N+1} \setminus S_j,
\end{equation}
and
\begin{equation}
\mu_w(S_i \cap S_j) = 0 \qquad \text{for all } i \ne j.
\end{equation}
Then
\begin{equation}
v_j \rightharpoonup 0 \qquad \text{weakly in } E.
\end{equation}
\end{lemma}

\begin{proof}
By hypothesis, there exists $M > 0$ such that
\begin{equation}
\|v_j\|_E \le M
\end{equation}
for all $j \in \mathbb{N}$.

For each $j$, since $v_j = 0$ $\mu_w$-almost everywhere on $\mathbb{H}^{N+1} \setminus S_j$, the locality principle for the intrinsic weighted weak gradient yields
\begin{equation}
\nabla v_j = 0 \quad \mu_w\text{-almost everywhere on } \mathbb{H}^{N+1} \setminus S_j.
\end{equation}
Hence both $v_j$ and $\nabla v_j$ vanish $\mu_w$-almost everywhere outside $S_j$.

Let $f \in E$ be arbitrary. Since $f \in H^1_{\mu_w}(\mathbb{H}^{N+1})$, one has
\begin{equation}
|f|^2 \in L^1_{\mu_w}(\mathbb{H}^{N+1})
\qquad \text{and} \qquad
|\nabla f|^2 \in L^1_{\mu_w}(\mathbb{H}^{N+1}).
\end{equation}
Define
\begin{equation}
a_j = \int_{S_j} \bigl(|f|^2 + |\nabla f|^2\bigr)\,d\mu_w.
\end{equation}
Each $a_j$ is a nonnegative real number. Since the sets $(S_j)$ are pairwise disjoint up to $\mu_w$-null sets, countable additivity of the Lebesgue integral gives
\begin{equation}
\sum_{j=1}^{\infty} a_j
=
\int_{\bigcup_{j=1}^{\infty} S_j} \bigl(|f|^2 + |\nabla f|^2\bigr)\,d\mu_w
\le
\int_{\mathbb{H}^{N+1}} \bigl(|f|^2 + |\nabla f|^2\bigr)\,d\mu_w
=
\|f\|_E^2.
\end{equation}
Since $\|f\|_E^2 < \infty$, the series $\sum_{j=1}^{\infty} a_j$ converges, and therefore
\begin{equation}
a_j \to 0
\qquad \text{as } j \to \infty.
\end{equation}

Because $v_j$ and $\nabla v_j$ vanish $\mu_w$-almost everywhere outside $S_j$, one has
\begin{equation}
\langle v_j, f \rangle_E
=
\int_{S_j} \bigl(v_j f + \nabla v_j \cdot \nabla f\bigr)\,d\mu_w.
\end{equation}
Applying the Cauchy-Schwarz inequality in the product space
\begin{equation}
L^2_{\mu_w}(S_j) \times \bigl(L^2_{\mu_w}(S_j)\bigr)^{N+1},
\end{equation}
we obtain
\begin{equation}
|\langle v_j, f \rangle_E|
\le
\left(
\int_{S_j} \bigl(|v_j|^2 + |\nabla v_j|^2\bigr)\,d\mu_w
\right)^{1/2}
\left(
\int_{S_j} \bigl(|f|^2 + |\nabla f|^2\bigr)\,d\mu_w
\right)^{1/2}.
\end{equation}
Thus
\begin{equation}
|\langle v_j, f \rangle_E|
\le
\|v_j\|_{E(S_j)} \sqrt{a_j}.
\end{equation}
Since $S_j \subset \mathbb{H}^{N+1}$, it follows that
\begin{equation}
\|v_j\|_{E(S_j)} \le \|v_j\|_E \le M.
\end{equation}
Therefore
\begin{equation}
0 \le |\langle v_j, f \rangle_E| \le M \sqrt{a_j}.
\end{equation}
As $a_j \to 0$, the squeeze theorem yields
\begin{equation}
\langle v_j, f \rangle_E \to 0
\qquad \text{as } j \to \infty.
\end{equation}

Since $f \in E$ was arbitrary, this is precisely the definition of weak convergence of $v_j$ to $0$ in the Hilbert space $E$. Therefore
\begin{equation}
v_j \rightharpoonup 0
\qquad \text{weakly in } E.
\end{equation}
\end{proof}

\section{Functional Inequalities and Local Tools}\label{sec:Functional Inequalities and Local Tools}

The global compactness property naturally depends upon local compactness coupled with uniform control at the boundaries of the domain. In Section~\ref{sec:Functional Inequalities and Local Tools}, we restrict our attention to bounded, Lipschitz-regular domains that are strictly separated from the singular boundary. It is shown that the strict local regularity of the weight guarantees local compactness via the classical Rellich-Kondrachov theorem. To handle the inevitable boundary degeneracy, a one-dimensional weighted Hardy inequality is derived , alongside scale-invariant annular Poincaré inequalities. These localized functional inequalities are indispensable; they provide the strict bounds required to control the energy of functions as $y \rightarrow 0$.

\begin{definition}[Lipschitz-regular]
Let $\mathbb{H}^{N+1}=\mathbb{R} \times \mathbb{R}_{+}^N$. Let $\mathfrak{O}$ be the collection of open subsets $\Omega \subset \mathbb{H}^{N+1}$. An open set $\Omega \in \mathfrak{O}$ is said to be Lipschitz-regular if for every point $p \in \partial \Omega$, there exists a neighborhood $U$ of $p$ and a bijective bi-Lipschitz mapping $F: U \rightarrow B(0,1)$ such that $F(U \cap \Omega)$ corresponds to the positive half-ball.
\end{definition}

\begin{definition}[Interior Approximable Weighted Sobolev Space]\label{def:Interior Approximable Weighted Sobolev Space}
Let $c \leq-1$. We define the space of infinitely differentiable functions with compact support strictly separated from the boundary $\Gamma=\{y=0\}$ as:

$$
C_c^{\infty}\left(\mathbb{H}^{N+1}\right):=\left\{v \in C^{\infty}\left(\overline{\mathbb{H}^{N+1}}\right): \operatorname{supp}(v) \cap \Gamma=\emptyset\right\}
$$

Let $H_{\mu_w}^1$ be the completion of the vector space of smooth functions under the norm:

$$
\|u\|_{H_{\mu_w}^1}:=\left(\int_{\mathbb{H}^{N+1}}\left(|u|^2+|\nabla u|^2\right) d \mu_w\right)^{1 / 2}
$$

We define $\mathcal{V}_0 \equiv H_{\mu_w, 0}^1\left(\mathbb{H}^{N+1}\right)$ as the topological closure of $C_c^{\infty}\left(\mathbb{H}^{N+1}\right)$ strictly within the topology induced by the norm of $H_{\mu_w}^1$. For any $u \in \mathcal{V}_0$, we say that $u$ is "interior approximable."
\end{definition}

\begin{remark}
The construction of the space $\mathcal{V}_0$ strictly bypasses the complex trace theory associated with the singular boundary $\Gamma = \{y=0\}$. By taking the closure of smooth functions whose supports are strictly separated from $\Gamma$, the definition isolates the functional analytic behavior of the domain's interior, providing a rigorous framework for functions that "vanish" at the singularity without requiring a classical boundary trace.
\end{remark}

\begin{proposition}[Strict local compactness]\label{prop:Strict local compactness}
Let $\mathbb{H}^{N+1} = \mathbb{R}^{N} \times (0,\infty)$. Let $\Omega \subset \mathbb{H}^{N+1}$ be a bounded Lipschitz open set such that $\operatorname{dist}(\Omega,\{y=0\})>0$. Let $\mu_w$ be the measure defined by $d\mu_w = y^{c}\phi(|z|)\,dz$. Then the canonical embedding
\begin{equation}
\iota_w : H_{\mu_w}^{1}(\Omega) \to L_{\mu_w}^{2}(\Omega)
\end{equation}
is compact.
\end{proposition}

\begin{proof}
Because $\Omega$ is bounded, there exists $R_{\max} > 0$ such that $|z| \le R_{\max}$ for all $z \in \Omega$. The separation condition implies the existence of $\delta > 0$ such that $y \ge \delta$ for all $z = (y,x) \in \Omega$. Since $w(z) = y^{c}\phi(|z|)$ and $\phi$ is locally bounded, there exist constants $\lambda_{\min}, \lambda_{\max} > 0$ such that
\begin{equation}
0 < \lambda_{\min} \le w(z) \le \lambda_{\max} < \infty
\quad \text{for almost every } z \in \Omega.
\end{equation}

Let $V_H$ be the vector space of weakly differentiable functions on $\Omega$. Define
\begin{equation}
E_1 = (V_H, \|\cdot\|_{H^{1}(\Omega)}), 
\qquad
E_2 = (V_H, \|\cdot\|_{H_{\mu_w}^{1}(\Omega)}).
\end{equation}
Similarly, let $V_L$ be the vector space of square-integrable functions on $\Omega$, and define
\begin{equation}
F_1 = (V_L, \|\cdot\|_{L^{2}(\Omega)}), 
\qquad
F_2 = (V_L, \|\cdot\|_{L_{\mu_w}^{2}(\Omega)}).
\end{equation}

From the bounds on $w$, for all $u \in V_H$,
\begin{equation}
\sqrt{\lambda_{\min}} \|u\|_{H^{1}(\Omega)} 
\le 
\|u\|_{H_{\mu_w}^{1}(\Omega)} 
\le 
\sqrt{\lambda_{\max}} \|u\|_{H^{1}(\Omega)}.
\end{equation}
Hence the identity map $J_H : E_2 \to E_1$ defined by $J_H(u) = u$ is a bounded linear bijection with bounded inverse.

Similarly, for all $u \in V_L$,
\begin{equation}
\sqrt{\lambda_{\min}} \|u\|_{L^{2}(\Omega)} 
\le 
\|u\|_{L_{\mu_w}^{2}(\Omega)} 
\le 
\sqrt{\lambda_{\max}} \|u\|_{L^{2}(\Omega)}.
\end{equation}
Thus the identity map $J_L : F_1 \to F_2$ defined by $J_L(u) = u$ is a bounded linear bijection with bounded inverse.

Let $\iota : E_1 \to F_1$ denote the classical Sobolev embedding. Since $\Omega$ is bounded with Lipschitz boundary, the Rellich--Kondrachov theorem implies that $\iota$ is compact.

Define $\iota_w : E_2 \to F_2$ by $\iota_w(u) = u$. Then
\begin{equation}
\iota_w = J_L \circ \iota \circ J_H.
\end{equation}
The operator $J_H$ is bounded from $E_2$ to $E_1$, the operator $\iota$ is compact from $E_1$ to $F_1$, and $J_L$ is bounded from $F_1$ to $F_2$. Therefore, the composition $J_L \circ \iota \circ J_H$ is compact. Hence $\iota_w$ is compact.
\end{proof}

\begin{lemma}[Local Compactness of Weighted Sobolev Embeddings]\label{lem:Local Compactness of Weighted Sobolev Embeddings}
Let $\Omega \subset \mathbb{H}^{N+1}$ be an open, bounded set. If the boundary $\partial \Omega$ is a manifold of class $C^\infty$, then $\Omega$ is Lipschitz-regular.

Let $(\Omega_m)_{m \in \mathbb{N}}$ be an admissible smooth exhaustion of $\mathbb{H}^{N+1}$. Let $\mu_w$ be a measure given by $d\mu_w = w(z)\,dz$. Assume $w(z)$ satisfies the locally regular condition. Then for every $m \in \mathbb{N}$, the canonical linear embedding
\begin{equation}
\iota_m : H_{\mu_w}^1(\Omega_m) \hookrightarrow L_{\mu_w}^2(\Omega_m)
\end{equation}
is a compact operator.
\end{lemma}

\begin{proof}
By the definition of a $C^\infty$ manifold boundary in $\mathbb{R}^{N+1}$, for every point $p \in \partial \Omega$, there exists a neighborhood $U$ of $p$ and a $C^\infty$ diffeomorphism $\Phi : U \to B(0,1)$ such that $\Phi(U \cap \Omega)$ corresponds to the positive half-ball. Since $\overline{U}$ may be chosen compact and $\Phi$ is of class $C^\infty$, both $\Phi$ and $\Phi^{-1}$ have uniformly bounded Jacobians on $U$. By the mean value inequality, $\Phi$ is a bijective bi-Lipschitz mapping. Hence $\Omega$ is Lipschitz-regular.

Fix $m \in \mathbb{N}$. By the properties of the admissible smooth exhaustion, $\overline{\Omega_m}$ is a closed and bounded subset of $\mathbb{R}^{N+1}$ and therefore compact. Moreover, there exists $\delta_m > 0$ such that $y \ge \delta_m$ for all $z \in \overline{\Omega_m}$, so $\overline{\Omega_m}$ is separated from the set $\{y = 0\}$.

The weight is given by $w(z) = y^c \phi(|z|)$. The functions $z \mapsto y^c$ and $z \mapsto |z|$ are continuous on the compact set $\overline{\Omega_m}$, hence bounded. By local regularity, $\phi(|z|)$ is strictly positive and bounded away from zero and infinity on compact intervals. Therefore there exist constants $0 < \lambda_m \le \Lambda_m < \infty$ such that
\begin{equation}
\lambda_m \le w(z) \le \Lambda_m
\quad \text{for almost every } z \in \Omega_m.
\end{equation}

Let $V_H$ be the vector space of weakly differentiable functions on $\Omega_m$. Consider the two Hilbert space structures
\begin{equation}
E_{\mathrm{unweighted}} = (V_H, \|\cdot\|_{H^1(\Omega_m)}), \quad
E_{\mathrm{weighted}} = (V_H, \|\cdot\|_{H_{\mu_w}^1(\Omega_m)}).
\end{equation}
From the bounds on $w$, for all $u \in V_H$ one has
\begin{equation}
\sqrt{\lambda_m}\,\|u\|_{H^1(\Omega_m)}
\le \|u\|_{H_{\mu_w}^1(\Omega_m)}
\le \sqrt{\Lambda_m}\,\|u\|_{H^1(\Omega_m)}.
\end{equation}
Define the identity map $I_H : E_{\mathrm{weighted}} \to E_{\mathrm{unweighted}}$ by $I_H(u) = u$. This map is linear, continuous, and bijective.

Similarly, let $V_L$ be the space of square-integrable functions on $\Omega_m$, and define
\begin{equation}
F_{\mathrm{unweighted}} = L^2(\Omega_m), \quad
F_{\mathrm{weighted}} = L_{\mu_w}^2(\Omega_m).
\end{equation}
The identity map $I_L : F_{\mathrm{unweighted}} \to F_{\mathrm{weighted}}$ given by $I_L(u) = u$ satisfies
\begin{equation}
\|I_L(u)\|_{F_{\mathrm{weighted}}}
\le \sqrt{\Lambda_m}\,\|u\|_{F_{\mathrm{unweighted}}},
\end{equation}
so $I_L$ is linear, continuous, and bijective.

Since $\partial \Omega_m$ is of class $C^\infty$, the domain $\Omega_m$ is Lipschitz-regular. The classical Rellich--Kondrachov theorem implies that the embedding
\begin{equation}
\iota_{RK} : E_{\mathrm{unweighted}} \hookrightarrow F_{\mathrm{unweighted}}
\end{equation}
is compact.

By Proposition~\ref{prop:Strict local compactness}, the weighted embedding $\iota_m : E_{\mathrm{weighted}} \hookrightarrow F_{\mathrm{weighted}}$ can be written as
\begin{equation}
\iota_m = I_L \circ \iota_{RK} \circ I_H.
\end{equation}
As the composition of a continuous operator, a compact operator, and a continuous operator, $\iota_m$ is compact.
\end{proof}

\begin{proposition}[1D Weighted Hardy Inequality]\label{prop:1D Weighted Hardy Inequality}
Let $c > -1$ with $c \neq 1$. For any $u \in C_{c}^{\infty}(\overline{\mathbb{H}^{N+1}})$ that vanishes in a neighborhood of $\{y=0\}$, the following inequality holds:
\begin{equation}
\int_{\mathbb{H}^{N+1}}\frac{u^{2}}{y^{2}}y^{c}\,dz \le \frac{4}{(c-1)^{2}}\int_{\mathbb{H}^{N+1}}|\partial_{y}u|^{2}y^{c}\,dz.
\end{equation}
\end{proposition}

\begin{proof}
Let $u \in C_{c}^{\infty}(\overline{\mathbb{H}^{N+1}})$ such that $u$ vanishes in a neighborhood of $\{y=0\}$. By compact support, for any fixed $x \in \mathbb{R}^{N}$ there exists a closed interval $[a_{x}, b_{x}] \subset (0,\infty)$ such that $u(y,x)=0$ for all $y \notin [a_{x}, b_{x}]$.

By Fubini's theorem,
\begin{equation}
\int_{\mathbb{H}^{N+1}} \frac{u^{2}}{y^{2}} y^{c}\,dz
=
\int_{\mathbb{R}^{N}} \left( \int_{0}^{\infty} u^{2} y^{c-2}\,dy \right) dx.
\end{equation}

Fix $x \in \mathbb{R}^{N}$. For $\lambda \in \mathbb{R}$, consider the non-negative quadratic form
\begin{equation}
0 \le \int_{0}^{\infty} \left( \partial_{y}u + \lambda \frac{u}{y} \right)^{2} y^{c}\,dy.
\end{equation}

Expanding,
\begin{equation}
0 \le \int_{0}^{\infty} (\partial_{y}u)^{2} y^{c}\,dy
+ 2\lambda \int_{0}^{\infty} u (\partial_{y}u) y^{c-1}\,dy
+ \lambda^{2} \int_{0}^{\infty} u^{2} y^{c-2}\,dy.
\end{equation}

Since $2u(\partial_{y}u)=\partial_{y}(u^{2})$,
\begin{equation}
2\lambda \int_{0}^{\infty} u (\partial_{y}u) y^{c-1}\,dy
=
\lambda \int_{0}^{\infty} \partial_{y}(u^{2}) y^{c-1}\,dy.
\end{equation}

Integration by parts gives
\begin{equation}
\lambda \int_{0}^{\infty} \partial_{y}(u^{2}) y^{c-1}\,dy
=
\lambda \left[ u^{2} y^{c-1} \right]_{0}^{\infty}
- \lambda(c-1) \int_{0}^{\infty} u^{2} y^{c-2}\,dy.
\end{equation}

Because $\operatorname{supp}(u(\cdot,x)) \subset [a_{x}, b_{x}]$ with $a_{x}>0$ and $b_{x}<\infty$, the boundary term vanishes, and hence
\begin{equation}
2\lambda \int_{0}^{\infty} u (\partial_{y}u) y^{c-1}\,dy
=
-\lambda(c-1) \int_{0}^{\infty} u^{2} y^{c-2}\,dy.
\end{equation}

Substituting,
\begin{equation}
0 \le \int_{0}^{\infty} (\partial_{y}u)^{2} y^{c}\,dy
+ (\lambda^{2} - \lambda(c-1)) \int_{0}^{\infty} u^{2} y^{c-2}\,dy.
\end{equation}

Minimizing the polynomial $\lambda^{2} - \lambda(c-1)$ gives $\lambda = \frac{c-1}{2}$, with minimum value $-\frac{(c-1)^{2}}{4}$. Therefore,
\begin{equation}
0 \le \int_{0}^{\infty} (\partial_{y}u)^{2} y^{c}\,dy
- \frac{(c-1)^{2}}{4} \int_{0}^{\infty} u^{2} y^{c-2}\,dy.
\end{equation}

Rearranging,
\begin{equation}
\frac{(c-1)^{2}}{4} \int_{0}^{\infty} u^{2} y^{c-2}\,dy
\le
\int_{0}^{\infty} (\partial_{y}u)^{2} y^{c}\,dy.
\end{equation}

Since $c \neq 1$, multiplying by $\frac{4}{(c-1)^{2}}$ yields
\begin{equation}
\int_{0}^{\infty} \frac{u^{2}}{y^{2}} y^{c}\,dy
\le
\frac{4}{(c-1)^{2}} \int_{0}^{\infty} (\partial_{y}u)^{2} y^{c}\,dy.
\end{equation}

Finally, integrating over $\mathbb{R}^{N}$ and applying Fubini-Tonelli completes the proof.
\end{proof}

\begin{proposition}[Scale-invariant annular Poincaré inequality]\label{prop:Scale-invariant annular Poincaré inequality}
Let $N \in \mathbb{N}$ and for $R > 0$ define
\begin{equation}
A_R = \{ z \in \mathbb{R}^{N+1} : R < |z| < 2R \}.
\end{equation}
Let $\lambda$ denote the Lebesgue measure. Then for every $u \in H^1(A_R)$,
\begin{equation}
\int_{A_R} |u - u_{A_R}|^2 \, dz \le C R^2 \int_{A_R} |\nabla u|^2 \, dz,
\end{equation}
where
\begin{equation}
u_{A_R} = \frac{1}{\lambda(A_R)} \int_{A_R} u \, dz,
\end{equation}
and $C > 0$ is a constant independent of $R$ and $u$.
\end{proposition}

\begin{proof}
Define the reference annulus
\begin{equation}
A_1 = \{ y \in \mathbb{R}^{N+1} : 1 < |y| < 2 \}.
\end{equation}
The boundary of $A_1$ consists of two smooth spheres, hence $A_1$ is a bounded Lipschitz domain. Therefore there exists a constant $C_1 > 0$ such that for all $v \in H^1(A_1)$,
\begin{equation}
\int_{A_1} |v - v_{A_1}|^2 \, dy \le C_1 \int_{A_1} |\nabla v|^2 \, dy.
\end{equation}

Define the mapping $T_R : A_1 \to A_R$ by $T_R(y) = R y$. Let $u \in H^1(A_R)$ and define $v = u \circ T_R$. Then $v \in H^1(A_1)$ and
\begin{equation}
\nabla v(y) = R (\nabla u)(R y).
\end{equation}

Since $\lambda(A_R) = R^{N+1} \lambda(A_1)$, we compute
\begin{equation}
u_{A_R}
=
\frac{1}{\lambda(A_R)} \int_{A_R} u(z) \, dz
=
\frac{1}{R^{N+1} \lambda(A_1)} \left( R^{N+1} \int_{A_1} u(Ry) \, dy \right)
=
\frac{1}{\lambda(A_1)} \int_{A_1} v(y) \, dy
=
v_{A_1}.
\end{equation}

Using the change of variables,
\begin{equation}
\int_{A_R} |u(z) - u_{A_R}|^2 \, dz
=
R^{N+1} \int_{A_1} |u(Ry) - u_{A_R}|^2 \, dy
=
R^{N+1} \int_{A_1} |v(y) - v_{A_1}|^2 \, dy.
\end{equation}

Applying the Poincaré inequality on $A_1$,
\begin{equation}
R^{N+1} \int_{A_1} |v - v_{A_1}|^2 \, dy
\le
R^{N+1} C_1 \int_{A_1} |\nabla v|^2 \, dy.
\end{equation}

Using the relation between gradients,
\begin{equation}
\int_{A_R} |\nabla u(z)|^2 \, dz
=
R^{N+1} \int_{A_1} |\nabla u(Ry)|^2 \, dy
=
R^{N+1} \int_{A_1} R^{-2} |\nabla v(y)|^2 \, dy
=
R^{N-1} \int_{A_1} |\nabla v(y)|^2 \, dy.
\end{equation}

Therefore,
\begin{equation}
R^{N+1} \int_{A_1} |\nabla v|^2 \, dy
=
R^2 \int_{A_R} |\nabla u|^2 \, dz.
\end{equation}

Combining the estimates,
\begin{equation}
\int_{A_R} |u - u_{A_R}|^2 \, dz \le C_1 R^2 \int_{A_R} |\nabla u|^2 \, dz.
\end{equation}

Setting $C = C_1$ completes the proof.
\end{proof}

\begin{lemma}[Scale-invariant annular Poincaré inequality]\label{lem:Scale-invariant annular Poincaré inequality}
Let $N \in \mathbb{N}$ and set $d = N+1$. For any $R > 0$, define
\begin{equation}
A_R = \{ z \in \mathbb{R}^{N+1} : R < |z| < 2R \}.
\end{equation}
Let $\lambda$ denote the Lebesgue measure. Then for every $u \in H^1(A_R)$,
\begin{equation}
\int_{A_R} |u - u_{A_R}|^2 \, dz \le C R^2 \int_{A_R} |\nabla u|^2 \, dz,
\end{equation}
where
\begin{equation}
u_{A_R} = \frac{1}{\lambda(A_R)} \int_{A_R} u \, dz,
\end{equation}
and $C > 0$ is a constant independent of $R$ and $u$.
\end{lemma}

\begin{proof}
Define the reference annulus
\begin{equation}
A_1 = \{ y \in \mathbb{R}^{N+1} : 1 < |y| < 2 \}.
\end{equation}
The boundary of $A_1$ consists of two smooth spheres, hence $A_1$ is a bounded Lipschitz domain. Therefore there exists a constant $C_1 > 0$ such that for all $v \in H^1(A_1)$,
\begin{equation}
\int_{A_1} |v - v_{A_1}|^2 \, dy \le C_1 \int_{A_1} |\nabla v|^2 \, dy.
\end{equation}

Define the homothety $T_R : A_1 \to A_R$ by $T_R(y) = R y$. Let $u \in H^1(A_R)$ and define $v = u \circ T_R$. Then $v \in H^1(A_1)$ and
\begin{equation}
\nabla v(y) = R (\nabla u)(R y).
\end{equation}

Since $\lambda(A_R) = R^d \lambda(A_1)$, we compute
\begin{equation}
u_{A_R}
=
\frac{1}{\lambda(A_R)} \int_{A_R} u(z) \, dz
=
\frac{1}{R^d \lambda(A_1)} \left( R^d \int_{A_1} u(Ry) \, dy \right)
=
\frac{1}{\lambda(A_1)} \int_{A_1} v(y) \, dy
=
v_{A_1}.
\end{equation}

Using the change of variables,
\begin{equation}
\int_{A_R} |u(z) - u_{A_R}|^2 \, dz
=
R^d \int_{A_1} |u(Ry) - u_{A_R}|^2 \, dy
=
R^d \int_{A_1} |v(y) - v_{A_1}|^2 \, dy.
\end{equation}

Applying the Poincaré inequality on $A_1$,
\begin{equation}
R^d \int_{A_1} |v - v_{A_1}|^2 \, dy
\le
R^d C_1 \int_{A_1} |\nabla v|^2 \, dy.
\end{equation}

Using the relation between gradients,
\begin{equation}
\int_{A_R} |\nabla u(z)|^2 \, dz
=
R^d \int_{A_1} |\nabla u(Ry)|^2 \, dy
=
R^d \int_{A_1} R^{-2} |\nabla v(y)|^2 \, dy
=
R^{d-2} \int_{A_1} |\nabla v(y)|^2 \, dy.
\end{equation}

Therefore,
\begin{equation}
R^d \int_{A_1} |\nabla v|^2 \, dy
=
R^2 \int_{A_R} |\nabla u|^2 \, dz.
\end{equation}

Combining the estimates,
\begin{equation}
\int_{A_R} |u - u_{A_R}|^2 \, dz \le C_1 R^2 \int_{A_R} |\nabla u|^2 \, dz.
\end{equation}

Setting $C = C_1$ completes the proof.
\end{proof}

\section{Tail Control and Decay Estimates}\label{sec:Tail Control and Decay Estimates}

Having secured local compactness, the primary analytic difficulty resides in the non-compactness of the spatial infinity. Section~\ref{sec:Tail Control and Decay Estimates} is devoted to the asymptotic behavior of functions in the weighted Sobolev space. By exploiting measurable Lyapunov potentials, we establish uniform tail mass control directly from the coercive tail inequality. The central result of this section is of a constructive nature: Lemma~\ref{lem:Explicit Lyapunov Construction} proves that if the Sobolev unit ball is tail tight, one can explicitly construct a corresponding Lyapunov potential that enforces this coercivity. This structural equivalence bridges the gap between abstract topological tightness and explicit potential theory.

\begin{lemma}[Tail control from coercivity]\label{lem:Tail control from coercivity}
Assume Definition~\ref{def:tail coercivity}. Then there exist $R_{0} \ge 1$ and a function $\Theta(R) \downarrow 0$ such that for all $R \ge R_{0}$ and all $u \in C_{c}^{\infty}(\mathbb{H}^{N+1})$,
\begin{equation}
\int_{\{ |z| > 2R \}} u^{2}\, d\mu_{w}
\ \le\
\Theta(R)\,
\int_{\{ |z| > R \}} \big( |\nabla u|^{2} + u^{2} \big)\, d\mu_{w}.
\end{equation}
\end{lemma}

\begin{proof}
By assumption, there exist a function $\Psi : [0,\infty) \to [0,\infty)$, with $\Psi(\rho) \uparrow \infty$ as $\rho \to \infty$, and a constant $C > 0$ such that for all $v \in C_{c}^{\infty}(\overline{\mathbb{H}^{N+1}})$
(vanishing near $y = 0$ if $c \le -1$),
\begin{equation}
\int_{\mathbb{H}^{N+1}} \Psi(|z|)\, v^{2}\, d\mu_{w}
\ \le\
C \int_{\mathbb{H}^{N+1}} \big( |\nabla v|^{2} + v^{2} \big)\, d\mu_{w}.
\end{equation}
Since $u \in C_{c}^{\infty}(\mathbb{H}^{N+1})$ vanishes near $y = 0$, it is admissible by definition~\ref{def:Interior Approximable Weighted Sobolev Space}.

Fix $R \ge 1$ and choose $\eta_{R} \in C_{c}^{\infty}(\mathbb{R}^{N+1})$ radial such that
\begin{equation}
0 \le \eta_{R} \le 1, \qquad
\eta_{R}(z) = 0 \text{ if } |z| \le R, \qquad
\eta_{R}(z) = 1 \text{ if } |z| \ge 2R, \qquad
|\nabla \eta_{R}(z)| \le \frac{C_{0}}{R} \text{ for } R < |z| < 2R.
\end{equation}

Set $v := \eta_{R}\, u$.  
Then $v \in C_{c}^{\infty}(\mathbb{H}^{N+1})$. It follows by assumption that

\begin{equation}
\int_{\mathbb{H}^{N+1}} \Psi(|z|)\, \eta_{R}^{2} u^{2}\, d\mu_{w}
\ \le\
C \int_{\mathbb{H}^{N+1}}
\big( |\nabla(\eta_{R} u)|^{2} + \eta_{R}^{2} u^{2} \big)\, d\mu_{w}.
\end{equation}

On $\{ |z| > 2R \}$, $\eta_{R} \equiv 1$ and $\Psi(|z|) \ge \Psi(2R)$, so
\begin{equation}
\Psi(2R) \int_{\{ |z| > 2R \}} u^{2}\, d\mu_{w}
\ \le\
\int_{\mathbb{H}^{N+1}} \Psi(|z|)\, \eta_{R}^{2} u^{2}\, d\mu_{w}.
\end{equation}

The product rule gives
\begin{equation}
|\nabla(\eta_{R} u)|^{2}
= \eta_{R}^{2} |\nabla u|^{2}
+ 2 \eta_{R} u\, \nabla \eta_{R} \!\cdot\! \nabla u
+ u^{2} |\nabla \eta_{R}|^{2}.
\end{equation}
For $\varepsilon \in (0,1]$, the mixed term satisfies
\begin{equation}
2 |\eta_{R} u\, \nabla \eta_{R} \!\cdot\! \nabla u|
\le
\varepsilon\, \eta_{R}^{2} |\nabla u|^{2}
+ \frac{1}{\varepsilon} u^{2} |\nabla \eta_{R}|^{2}.
\end{equation}
Thus
\begin{equation}
|\nabla(\eta_{R} u)|^{2}
\le
(1 + \varepsilon)\, \eta_{R}^{2} |\nabla u|^{2}
+ \Big( 1 + \frac{1}{\varepsilon} \Big) u^{2} |\nabla \eta_{R}|^{2}.
\end{equation}

The gradient satisfies $\operatorname{supp} \nabla \eta_{R} \subset A_{R} := \{ R < |z| < 2R \}$ and
$|\nabla \eta_{R}|^{2} \le C_{0}^{2} / R^{2}$ there.

Combining the previous estimates,
\begin{equation}
\int_{\mathbb{H}^{N+1}}
\big( |\nabla(\eta_{R} u)|^{2} + \eta_{R}^{2} u^{2} \big)\, d\mu_{w}
\le
(1 + \varepsilon) \int_{\{ |z| > R \}} \eta_{R}^{2} |\nabla u|^{2}\, d\mu_{w}
+ \Big( 1 + \frac{1}{\varepsilon} \Big) \frac{C_{0}^{2}}{R^{2}} \int_{A_{R}} u^{2}\, d\mu_{w}
+ \int_{\{ |z| > R \}} \eta_{R}^{2} u^{2}\, d\mu_{w}.
\end{equation}

Since $0 \le \eta_{R} \le 1$ and $\eta_{R} = 0$ for $|z| \le R$,
\begin{equation}
\int_{\{ |z| > R \}} \eta_{R}^{2} |\nabla u|^{2}
\le \int_{\{ |z| > R \}} |\nabla u|^{2}, 
\qquad
\int_{\{ |z| > R \}} \eta_{R}^{2} u^{2}
\le \int_{\{ |z| > R \}} u^{2}.
\end{equation}

\begin{equation}
\Psi(2R) \int_{\{ |z| > 2R \}} u^{2}\, d\mu_{w}
\le
C \Big[
(1 + \varepsilon) \int_{\{ |z| > R \}} |\nabla u|^{2}\, d\mu_{w}
+ \int_{\{ |z| > R \}} u^{2}\, d\mu_{w}
+ \Big( 1 + \frac{1}{\varepsilon} \Big) \frac{C_{0}^{2}}{R^{2}} \int_{A_{R}} u^{2}\, d\mu_{w}
\Big].
\end{equation}

Since $A_{R} \subset \{ |z| > R \}$,
\begin{equation}
\int_{A_{R}} u^{2} \le \int_{\{ |z| > R \}} u^{2}.
\end{equation}
Choosing $\varepsilon = 1$, we obtain
\begin{equation}
\Psi(2R) \int_{\{ |z| > 2R \}} u^{2}\, d\mu_{w}
\le
C_{1} \int_{\{ |z| > R \}} \big( |\nabla u|^{2} + u^{2} \big)\, d\mu_{w}
+ \frac{C_{2}}{R^{2}} \int_{\{ |z| > R \}} u^{2}\, d\mu_{w}.
\end{equation}

Rewriting the previous inequality,
\begin{equation}
\int_{\{ |z| > 2R \}} u^{2}\, d\mu_{w}
\le
\underbrace{\frac{C_{1}}{\Psi(2R)}}_{=: \Theta_{1}(R)}
\int_{\{ |z| > R \}} (|\nabla u|^{2} + u^{2})\, d\mu_{w}
+
\underbrace{\frac{C_{2}}{R^{2} \Psi(2R)}}_{=: \Theta_{2}(R)}
\int_{\{ |z| > R \}} u^{2}\, d\mu_{w}.
\end{equation}

Since $\Psi(2R) \uparrow \infty$ as $R \to \infty$, both $\Theta_{1}(R)$ and $\Theta_{2}(R)$ decrease to $0$.  
Define $\Theta(R) := \Theta_{1}(R) + \Theta_{2}(R)$, which satisfies $\Theta(R) \downarrow 0$ as $R \to \infty$.

Select $R_{0} \ge 1$ large enough so that $\Psi(2R) > 0$ and
$\max\{ \Theta_{1}(R), \Theta_{2}(R) \} \le 1$ for all $R \ge R_{0}$.

Hence there exist $R_{0} \ge 1$ and $\Theta(R) \downarrow 0$ such that
\begin{equation}
\int_{\{ |z| > 2R \}} u^{2}\, d\mu_{w}
\le
\Theta(R)
\int_{\{ |z| > R \}} \big( |\nabla u|^{2} + u^{2} \big)\, d\mu_{w},
\end{equation}
for all $R \ge R_{0}$ and all $u \in C_{c}^{\infty}(\mathbb{H}^{N+1})$.
\end{proof}

\begin{remark}
The significance of this lemma lies in its localization of the abstract coercive property. By employing a standard radial cut-off argument, it explicitly translates the global structural coercivity hypothesis into a localized decay estimate at spatial infinity. The bounding function $\Theta(R)$ serves as a quantifiable rate of mass decay, strictly binding the $L_{\mu_w}^2$-mass in the far tail region $\{|z| > 2R\}$ to the Sobolev energy in the annular exterior $\{|z| > R\}$.
\end{remark}

\begin{lemma}[Lower Semi-Continuity of the Lyapunov Functional]\label{lem:Lower Semi-Continuity of the Lyapunov Functional}
Let $(\mathbb{H}^{N+1}, \mu_{w})$ be the measure space defined by $\mu_{w} = y^{c}\phi(|z|)\,dz$. Let $\Psi: [0,\infty) \to [0,\infty)$ be a Borel measurable function. If the coercive tail inequality
\begin{equation}
\int_{\mathbb{H}^{N+1}} \Psi(|z|) v^{2}\, d\mu_{w} \le C_{TC} \|v\|_{H_{\mu_{w}}^{1}}^{2}
\end{equation}
holds for all $v \in C_{c}^{\infty}(\mathbb{H}^{N+1})$ (vanishing near $y=0$ if $c \le -1$), then it holds for all $u \in H_{\mu_{w},0}^{1}(\mathbb{H}^{N+1})$.
\end{lemma}

\begin{proof}
Let $u \in H_{\mu_{w},0}^{1}(\mathbb{H}^{N+1})$. By the definition of closure, there exists a sequence $(v_{n})_{n \in \mathbb{N}} \subset C_{c}^{\infty}(\mathbb{H}^{N+1})$ such that $v_{n} \to u$ strongly in $H_{\mu_{w}}^{1}$.

Since the $H_{\mu_{w}}^{1}$ norm dominates the $L_{\mu_{w}}^{2}$ norm,
\begin{equation}
\lim_{n \to \infty} \|v_{n} - u\|_{L_{\mu_{w}}^{2}} = 0.
\end{equation}

Thus $(v_{n})$ is Cauchy in $L_{\mu_{w}}^{2}(\mathbb{H}^{N+1})$. By the Riesz-Fischer theorem, there exists a subsequence $(v_{n_{k}})$ such that
\begin{equation}
\lim_{k \to \infty} v_{n_{k}}(z) = u(z) \quad \text{for } \mu_{w}\text{-a.e. } z \in \mathbb{H}^{N+1}.
\end{equation}

Define
\begin{equation}
f_{k}(z) := \Psi(|z|)\, |v_{n_{k}}(z)|^{2}.
\end{equation}

Since $\Psi$ is Borel measurable and non-negative, and $v_{n_{k}}$ are measurable, each $f_{k}$ is non-negative and $\mu_{w}$-measurable. By pointwise convergence,
\begin{equation}
\lim_{k \to \infty} f_{k}(z) = \Psi(|z|)\, |u(z)|^{2} \quad \text{for } \mu_{w}\text{-a.e. } z.
\end{equation}

Applying Fatou's lemma,
\begin{equation}
\int_{\mathbb{H}^{N+1}} \Psi(|z|)\, |u|^{2}\, d\mu_{w}
\le
\liminf_{k \to \infty} \int_{\mathbb{H}^{N+1}} \Psi(|z|)\, |v_{n_{k}}|^{2}\, d\mu_{w}.
\end{equation}

By the hypothesis applied to $v_{n_{k}}$,
\begin{equation}
\int_{\mathbb{H}^{N+1}} \Psi(|z|)\, |v_{n_{k}}|^{2}\, d\mu_{w}
\le
C_{TC} \|v_{n_{k}}\|_{H_{\mu_{w}}^{1}}^{2}.
\end{equation}

Hence,
\begin{equation}
\liminf_{k \to \infty} \int_{\mathbb{H}^{N+1}} \Psi(|z|)\, |v_{n_{k}}|^{2}\, d\mu_{w}
\le
\liminf_{k \to \infty} C_{TC} \|v_{n_{k}}\|_{H_{\mu_{w}}^{1}}^{2}.
\end{equation}

Since the norm is continuous on $H_{\mu_{w}}^{1}$,
\begin{equation}
\lim_{k \to \infty} \|v_{n_{k}}\|_{H_{\mu_{w}}^{1}}^{2} = \|u\|_{H_{\mu_{w}}^{1}}^{2}.
\end{equation}

Therefore,
\begin{equation}
\liminf_{k \to \infty} C_{TC} \|v_{n_{k}}\|_{H_{\mu_{w}}^{1}}^{2}
=
C_{TC} \|u\|_{H_{\mu_{w}}^{1}}^{2}.
\end{equation}

Combining the inequalities,
\begin{equation}
\int_{\mathbb{H}^{N+1}} \Psi(|z|)\, |u|^{2}\, d\mu_{w}
\le
C_{TC} \|u\|_{H_{\mu_{w}}^{1}}^{2}.
\end{equation}

This proves the claim.
\end{proof}

\begin{remark}
This lemma provides an indispensable technical bridge. While the tail coercivity condition is naturally hypothesized on the dense subspace of smooth, compactly supported test functions, the subsequent abstract topological arguments require these bounds to hold strictly on the full completion space $H_{\mu_w,0}^1$. The proof elegantly exploits Fatou's lemma alongside the almost-everywhere pointwise convergence of a strongly Cauchy sequence, rigorously confirming that the Lyapunov potential behaves as a proper, weakly lower semi-continuous functional.
\end{remark}

\begin{lemma}[Tail Control via Lyapunov Coercivity]\label{lem:Tail Control via Lyapunov Coercivity}
Let $(\mathbb{H}^{N+1}, \mu_w)$ be the measurable space defined by $\mu_w = y^{c}\phi(|z|)\,dz$. Assume the weight admits a Borel measurable Lyapunov potential $\Psi: [0,\infty) \to [0,\infty)$ satisfying the Tail Coercivity condition. For any $R \ge 1$ sufficiently large such that the essential infimum
\begin{equation}
m_R := \operatorname{ess\,inf}_{|z|\ge R} \Psi(|z|) > 0,
\end{equation}
the following estimate holds for all $u \in H_{\mu_w}^{1}(\mathbb{H}^{N+1})$:
\begin{equation}
\int_{\{|z|>R\}} |u|^{2}\, d\mu_w \le \frac{C_{TC}}{m_R} \|u\|_{H_{\mu_w}^{1}}^{2}.
\end{equation}
\end{lemma}

\begin{proof}
Let $u \in H_{\mu_w}^{1}(\mathbb{H}^{N+1})$. By Lemma~\ref{lem:Lower Semi-Continuity of the Lyapunov Functional}, the coercive tail inequality extends to $u$:
\begin{equation}
\int_{\mathbb{H}^{N+1}} \Psi(|z|)\, |u|^{2}\, d\mu_w \le C_{TC} \|u\|_{H_{\mu_w}^{1}}^{2}.
\end{equation}

Define the measurable set
\begin{equation}
E_{\mathrm{out}} := \{ z \in \mathbb{H}^{N+1} : |z| > R \}.
\end{equation}

By definition of the essential infimum, for $\mu_w$-almost every $z \in E_{\mathrm{out}}$,
\begin{equation}
\Psi(|z|)\, |u(z)|^{2} \ge m_R\, |u(z)|^{2}.
\end{equation}

Integrating over $E_{\mathrm{out}}$ and using monotonicity of the Lebesgue integral,
\begin{equation}
\int_{E_{\mathrm{out}}} \Psi(|z|)\, |u|^{2}\, d\mu_w \ge m_R \int_{E_{\mathrm{out}}} |u|^{2}\, d\mu_w.
\end{equation}

Partitioning $\mathbb{H}^{N+1}$ into $E_{\mathrm{out}}$ and its complement,
\begin{equation}
\int_{\mathbb{H}^{N+1}} \Psi(|z|)\, |u|^{2}\, d\mu_w
=
\int_{E_{\mathrm{out}}} \Psi(|z|)\, |u|^{2}\, d\mu_w
+
\int_{\mathbb{H}^{N+1}\setminus E_{\mathrm{out}}} \Psi(|z|)\, |u|^{2}\, d\mu_w.
\end{equation}

Since the integrand is non-negative,
\begin{equation}
\int_{\mathbb{H}^{N+1}} \Psi(|z|)\, |u|^{2}\, d\mu_w
\ge
\int_{E_{\mathrm{out}}} \Psi(|z|)\, |u|^{2}\, d\mu_w.
\end{equation}

Combining the inequalities,
\begin{equation}
m_R \int_{E_{\mathrm{out}}} |u|^{2}\, d\mu_w
\le
\int_{E_{\mathrm{out}}} \Psi(|z|)\, |u|^{2}\, d\mu_w
\le
\int_{\mathbb{H}^{N+1}} \Psi(|z|)\, |u|^{2}\, d\mu_w
\le
C_{TC} \|u\|_{H_{\mu_w}^{1}}^{2}.
\end{equation}

Since $m_R > 0$, dividing yields
\begin{equation}
\int_{\{|z|>R\}} |u|^{2}\, d\mu_w \le \frac{C_{TC}}{m_R} \|u\|_{H_{\mu_w}^{1}}^{2}.
\end{equation}

This completes the proof.
\end{proof}

\begin{remark}
Here, the abstract mechanism of the Lyapunov potential is rendered entirely explicit. Provided the potential $\Psi$ genuinely diverges at spatial infinity, its strictly positive essential infimum $m_R$ on the exterior region acts as an effective confining parameter. The resulting estimate neatly isolates the tail mass strictly in terms of the total Sobolev norm, effectively proving that the presence of a coercive potential forces uniform asymptotic mass-vanishing for every element in the Sobolev space.
\end{remark}

\begin{definition}[tail tightness modulus]
Let $E=H_{\mu_w}^1\left(\mathbb{H}^{N+1}\right)$ and let $\mathcal{B}_E$ denote its closed unit ball. For any $R>0$, let $A_R=\left\{z \in \mathbb{H}^{N+1}:|z|>R\right\}$. The tail tightness modulus is the mapping $\omega:(0, \infty) \rightarrow \mathbb{R}^{+}$defined by:

$$
\omega(R):=\sup _{u \in \mathcal{B}_E} \int_{A_R}|u(z)|^2 d \mu_w
$$
\end{definition}

\begin{remark}
The introduction of the modulus $\omega(R)$ provides a rigorous analytic quantification of the topological tightness property. Instead of merely asserting that mass vanishes at infinity, this supremum defines the absolute worst-case mass decay profile for the closed unit ball $\mathcal{B}_E$. It translates an abstract geometric condition into a precise metric function, establishing the exact rate at which mass escapes the system as $R \rightarrow \infty$.
\end{remark}

\begin{lemma}[Explicit Lyapunov Construction]\label{lem:Explicit Lyapunov Construction}
Let $E = H_{\mu_w}^1(\mathbb{H}^{N+1})$ be equipped with the canonical norm $\|\cdot\|_E$, and let $\mathcal{B}_E$ denote its closed unit ball. Assume $\mathcal{B}_E$ is tail tight, meaning that the tail modulus
\begin{equation}
\omega(R) := \sup_{v \in \mathcal{B}_E} \int_{|z| > R} |v|^2 \, d\mu_w
\end{equation}
satisfies $\lim_{R \to \infty} \omega(R) = 0$. Then there exists a Borel measurable function $\Psi : [0,\infty) \to [0,\infty)$ such that $\Psi(\rho)$ is non-decreasing and $\Psi(\rho) \to \infty$ as $\rho \to \infty$, and a constant $C_{TC} > 0$ such that for all $u \in E$,
\begin{equation}
\int_{\mathbb{H}^{N+1}} \Psi(|z|)\, |u|^2 \, d\mu_w \le C_{TC} \|u\|_E^2.
\end{equation}
\end{lemma}

\begin{proof}
Since $\lim_{R \to \infty} \omega(R) = 0$, for every $\epsilon > 0$ there exists $R^* > 0$ such that for all $R \ge R^*$ one has $\omega(R) \le \epsilon$.

We construct a sequence $(R_k)_{k \ge 0}$ inductively. Set $R_0 = 0$. Given $R_k$, choose $R^*_{k+1} > 0$ such that $\omega(R^*_{k+1}) \le 4^{-(k+1)}$, and define
\begin{equation}
R_{k+1} := \max(R_k + 1, R^*_{k+1}).
\end{equation}
Then $R_{k+1} > R_k$ for all $k$, and $R_k \ge k$, hence $R_k \to \infty$. Since $\omega$ is non-increasing, it follows that $\omega(R_{k+1}) \le 4^{-(k+1)}$.

Define $\Psi : [0,\infty) \to [0,\infty)$ by
\begin{equation}
\Psi(\rho) := \sum_{k=1}^{\infty} 2^k \chi_{[R_k, R_{k+1})}(\rho).
\end{equation}
This function is Borel measurable as a countable sum of simple functions. If $\rho_1 \le \rho_2$, and $\rho_1 \in [R_m, R_{m+1})$, $\rho_2 \in [R_n, R_{n+1})$, then $n \ge m$, hence $\Psi(\rho_2) \ge \Psi(\rho_1)$. Thus $\Psi$ is non-decreasing. Moreover, for any $M > 0$, there exists $K$ such that $2^K > M$, and for all $\rho \ge R_K$ one has $\Psi(\rho) \ge 2^K > M$, so $\Psi(\rho) \to \infty$ as $\rho \to \infty$.

Let $v \in \mathcal{B}_E$. Then
\begin{equation}
\int_{\mathbb{H}^{N+1}} \Psi(|z|)\, |v|^2 \, d\mu_w
= \sum_{k=1}^{\infty} 2^k \int_{\{R_k \le |z| < R_{k+1}\}} |v|^2 \, d\mu_w.
\end{equation}
Since $\{R_k \le |z| < R_{k+1}\} \subset \{|z| > R_k\}$, we obtain
\begin{equation}
\int_{\{R_k \le |z| < R_{k+1}\}} |v|^2 \, d\mu_w
\le \int_{\{|z| > R_k\}} |v|^2 \, d\mu_w
\le \omega(R_k)
\le 4^{-k}.
\end{equation}
Therefore
\begin{equation}
\int_{\mathbb{H}^{N+1}} \Psi(|z|)\, |v|^2 \, d\mu_w
\le \sum_{k=1}^{\infty} 2^k \cdot 4^{-k}
= \sum_{k=1}^{\infty} 2^{-k}
= 1.
\end{equation}

Let $u \in E$. If $u = 0$, then the inequality holds trivially. If $u \neq 0$, define $v := \|u\|_E^{-1} u$, so that $\|v\|_E = 1$ and $v \in \mathcal{B}_E$. Then
\begin{equation}
\int_{\mathbb{H}^{N+1}} \Psi(|z|)\, \left( \frac{|u|}{\|u\|_E} \right)^2 \, d\mu_w \le 1.
\end{equation}
By homogeneity of the integral,
\begin{equation}
\frac{1}{\|u\|_E^2}
\int_{\mathbb{H}^{N+1}} \Psi(|z|)\, |u|^2 \, d\mu_w \le 1,
\end{equation}
and multiplying by $\|u\|_E^2$ yields
\begin{equation}
\int_{\mathbb{H}^{N+1}} \Psi(|z|)\, |u|^2 \, d\mu_w \le \|u\|_E^2.
\end{equation}

Thus the conclusion holds with $C_{TC} = 1$.
\end{proof}

\begin{remark}
This lemma constitutes the structural culmination of the section. It establishes a profound converse to the preceding decay estimates: the topological condition of tail tightness is not merely a consequence of coercivity, but it constructively guarantees the existence of a coercive Lyapunov potential. The explicit summation construction over annular shells directly maps the metric decay prescribed by $\omega(R)$ into a divergent, non-decreasing potential function $\Psi$. This completes a strict equivalence, elegantly reducing abstract compactness to verifiable potential theory.
\end{remark}

\section{Abstract Compactness and Regularity}\label{sec:Abstract Compactness and Regularity}

We now translate the established analytic estimates into topological properties. Section~\ref{sec:Abstract Compactness and Regularity} formulates the abstract criteria for compactness utilizing a topological exhaustion method via localized orthogonal projections. It is established that the compactness of the continuous linear embedding $\iota: E \rightarrow F$ is entirely equivalent to localized compactness combined with the uniform vanishing of the projection remainder on the unit ball. Furthermore, we prove the global density of smooth, compactly supported functions within the weighted space, which justifies the limiting arguments required for the subsequent global theorems.

\begin{lemma}[Precompactness criterion via projection exhaustion]\label{lem:Precompactness criterion via projection exhaustion}
Let $E$ be a Banach space with closed unit ball $B_E$. Let $F$ be a Hilbert space. Let $\iota : E \to F$ be a continuous linear embedding. Let $(P_m)_{m \in \mathbb{N}}$ be a sequence of orthogonal projections on $F$, and define $P_m^\perp = I - P_m$.

Assume that for every $m \in \mathbb{N}$, the operator $P_m \circ \iota : E \to F$ is compact, and that
\begin{equation}
\lim_{m \to \infty} \sup_{u \in B_E} \|P_m^\perp(\iota(u))\|_F = 0.
\end{equation}

Then the set $\iota(B_E)$ is totally bounded in $F$. In particular, $\iota$ is a compact operator.
\end{lemma}

\begin{proof}
Let $\epsilon > 0$. We construct a finite $\epsilon$-net for $\iota(B_E)$ in $F$.

By the convergence assumption, there exists $M \in \mathbb{N}$ such that for all $u \in B_E$,
\begin{equation}
\|P_M^\perp(\iota(u))\|_F < \frac{\epsilon}{2}.
\end{equation}

Since $P_M \circ \iota$ is compact, the set
\begin{equation}
S_M = (P_M \circ \iota)(B_E)
\end{equation}
is relatively compact in $F$. Because $F$ is complete, this implies that $S_M$ is totally bounded. Therefore, there exist finitely many elements $v_1, \dots, v_k \in F$ such that for every $u \in B_E$, there exists $j \in \{1,\dots,k\}$ satisfying
\begin{equation}
\|P_M(\iota(u)) - v_j\|_F < \frac{\epsilon}{2}.
\end{equation}

Let $u \in B_E$. Using the identity $\iota(u) = P_M(\iota(u)) + P_M^\perp(\iota(u))$, we obtain
\begin{equation}
\|\iota(u) - v_j\|_F
=
\|P_M(\iota(u)) + P_M^\perp(\iota(u)) - v_j\|_F.
\end{equation}

By the triangle inequality,
\begin{equation}
\|\iota(u) - v_j\|_F
\le
\|P_M(\iota(u)) - v_j\|_F + \|P_M^\perp(\iota(u))\|_F.
\end{equation}

Combining the estimates,
\begin{equation}
\|\iota(u) - v_j\|_F < \frac{\epsilon}{2} + \frac{\epsilon}{2} = \epsilon.
\end{equation}

Thus $\{v_1,\dots,v_k\}$ is a finite $\epsilon$-net for $\iota(B_E)$. Hence $\iota(B_E)$ is totally bounded in $F$. Since $F$ is complete, its closure is compact, and therefore $\iota$ is a compact operator.
\end{proof}

\begin{theorem}[Topological characterization of global compactness]\label{thm:Topological characterization of global compactness}
Let $E = H_{\mu_w}^1(\mathbb{H}^{N+1})$ and $F = L_{\mu_w}^2(\mathbb{H}^{N+1})$. Let $\iota : E \to F$ be the canonical linear embedding. Let $(K_m)_{m \in \mathbb{N}}$ be an exhaustion of $\mathbb{H}^{N+1}$ by compact sets. For each $m \in \mathbb{N}$, define $P_m : F \to F$ as the orthogonal projection given by multiplication by $\chi_{K_m}$, and set $P_m^\perp = I - P_m$.

Then $\iota$ is compact if and only if the following two conditions hold:

For every $m \in \mathbb{N}$, the operator $P_m \circ \iota : E \to F$ is compact.

\begin{equation}
\lim_{m \to \infty} \sup_{u \in B_E} \|P_m^\perp(\iota(u))\|_F = 0,
\end{equation}
where $B_E$ denotes the closed unit ball of $E$.
\end{theorem}

\begin{proof}
Assume first that for every $m \in \mathbb{N}$ the operator $P_m \circ \iota$ is compact and that
\begin{equation}
\lim_{m \to \infty} \sup_{u \in B_E} \|P_m^\perp(\iota(u))\|_F = 0.
\end{equation}
By Lemma~\ref{lem:Precompactness criterion via projection exhaustion}, the set $\iota(B_E)$ is totally bounded in $F$. Since $F$ is complete, the closure of $\iota(B_E)$ is compact, and therefore $\iota$ is a compact operator.

Conversely, assume that $\iota$ is compact. We show that both conditions hold.

Fix $m \in \mathbb{N}$. The operator $P_m$ is a continuous linear operator on $F$ with operator norm at most $1$. The composition of a continuous linear operator with a compact operator is compact, hence $P_m \circ \iota$ is compact.

Let $S = \overline{\iota(B_E)}$, which is a compact subset of $F$. Fix $\epsilon > 0$. Since $S$ is compact, it is totally bounded, and there exists a finite set $\{v_1, \dots, v_k\} \subset F$ such that for every $f \in S$ there exists $j \in \{1,\dots,k\}$ with
\begin{equation}
\|f - v_j\|_F < \frac{\epsilon}{2}.
\end{equation}

For each $j$, since $v_j \in F$, the function $|v_j|^2$ is integrable with respect to $\mu_w$. By absolute continuity of the integral, as $m \to \infty$,
\begin{equation}
\|P_m^\perp v_j\|_F \to 0.
\end{equation}
Since the set $\{v_1,\dots,v_k\}$ is finite, there exists $M \in \mathbb{N}$ such that for all $m \ge M$ and all $j$,
\begin{equation}
\|P_m^\perp v_j\|_F < \frac{\epsilon}{2}.
\end{equation}

Let $u \in B_E$. Then $\iota(u) \in S$, so there exists $j$ such that
\begin{equation}
\|\iota(u) - v_j\|_F < \frac{\epsilon}{2}.
\end{equation}
Using that $\|P_m^\perp\| \le 1$,
\begin{equation}
\|P_m^\perp(\iota(u))\|_F
=
\|P_m^\perp(\iota(u) - v_j + v_j)\|_F
\le
\|P_m^\perp(\iota(u) - v_j)\|_F + \|P_m^\perp v_j\|_F.
\end{equation}
Hence,
\begin{equation}
\|P_m^\perp(\iota(u))\|_F
\le
\|\iota(u) - v_j\|_F + \|P_m^\perp v_j\|_F
<
\frac{\epsilon}{2} + \frac{\epsilon}{2}
= \epsilon.
\end{equation}

Since $M$ is independent of $u$, it follows that
\begin{equation}
\sup_{u \in B_E} \|P_m^\perp(\iota(u))\|_F \le \epsilon
\quad \text{for all } m \ge M.
\end{equation}

Thus,
\begin{equation}
\lim_{m \to \infty} \sup_{u \in B_E} \|P_m^\perp(\iota(u))\|_F = 0.
\end{equation}

This completes the proof.
\end{proof}

\begin{proposition}[Global density for integrable weights]\label{prop:Global density for integrable weights}
Let the parameter satisfy $c > -1$. Then the closure of the space $C_c^\infty(\mathbb{H}^{N+1})$ in the topology of $H_{\mu_w}^1$ is equal to the entire space $H_{\mu_w}^1(\mathbb{H}^{N+1})$.
\end{proposition}

\begin{proof}
Let $u \in H_{\mu_w}^1(\mathbb{H}^{N+1})$. We construct a sequence in $C_c^\infty(\mathbb{H}^{N+1})$ converging strongly to $u$.

Fix $R > 1$. Let $\zeta_R \in C_c^\infty(\mathbb{R}^{N+1})$ be a radial cut-off function such that $0 \le \zeta_R \le 1$, $\zeta_R \equiv 1$ on the ball $B_R$, $\zeta_R \equiv 0$ on $\mathbb{R}^{N+1} \setminus B_{2R}$, and $|\nabla \zeta_R| \le R^{-1}$. Define
\begin{equation}
u^{(R)} = \zeta_R u.
\end{equation}
Since $u \in H_{\mu_w}^1$, the functions $|u|^2$ and $|\nabla u|^2$ are integrable with respect to $\mu_w$. By dominated convergence,
\begin{equation}
\lim_{R \to \infty} \|u^{(R)} - u\|_{H_{\mu_w}^1}^2
=
\lim_{R \to \infty} \int_{\{|z|>R\}} \left(|u|^2 + |\nabla u|^2\right)\, d\mu_w
= 0.
\end{equation}

Fix $R$. For $\epsilon \in (0,1)$, choose $\eta_\epsilon \in C^\infty([0,\infty))$ such that $0 \le \eta_\epsilon \le 1$, $\eta_\epsilon \equiv 0$ on $[0,\epsilon]$, $\eta_\epsilon \equiv 1$ on $[2\epsilon,\infty)$, and $|\eta_\epsilon'| \le \epsilon^{-1}$. Define
\begin{equation}
u^{(R,\epsilon)} = \eta_\epsilon u^{(R)}.
\end{equation}
Then
\begin{equation}
u^{(R)} - u^{(R,\epsilon)} = (1 - \eta_\epsilon) u^{(R)},
\end{equation}
and this difference is supported in the strip $\{0 < y < 2\epsilon\}$.

Since $c > -1$, the measure $\mu_w$ is locally finite near the boundary. By absolute continuity of the integral,
\begin{equation}
\lim_{\epsilon \to 0} \int_{0<y<2\epsilon} |u^{(R)}|^2\, d\mu_w = 0.
\end{equation}

For the gradient, using the product rule,
\begin{equation}
\nabla\bigl((1 - \eta_\epsilon) u^{(R)}\bigr)
=
(1 - \eta_\epsilon)\nabla u^{(R)} - u^{(R)} \eta_\epsilon' e_y.
\end{equation}
The first term converges to zero by absolute continuity. For the second term,
\begin{equation}
\int_{\mathbb{H}^{N+1}} |u^{(R)} \eta_\epsilon'|^2\, d\mu_w
\le
\frac{1}{\epsilon^2} \int_{\{0<y<2\epsilon\}} |u^{(R)}|^2 y^c\, dz
\le
4 \int_{\{0<y<2\epsilon\}} \frac{|u^{(R)}|^2}{y^2} y^c\, dz.
\end{equation}
Applying the one-dimensional weighted Hardy inequality (Proposition~\ref{prop:Scale-invariant annular Poincaré inequality} to $u^{(R)}$,
\begin{equation}
\int_{\{0<y<2\epsilon\}} \frac{|u^{(R)}|^2}{y^2} y^c\, dz
\le
C \int_{\{0<y<2\epsilon\}} \left(|\partial_y u^{(R)}|^2 + |u^{(R)}|^2\right) y^c\, dz.
\end{equation}
As $\epsilon \to 0$, the right-hand side converges to zero by absolute continuity. Hence
\begin{equation}
\lim_{\epsilon \to 0} \|u^{(R,\epsilon)} - u^{(R)}\|_{H_{\mu_w}^1} = 0.
\end{equation}

The function $u^{(R,\epsilon)}$ has compact support contained in a set strictly separated from $\{y=0\}$. Let $\Omega \subset \mathbb{H}^{N+1}$ be an open bounded set containing this support with $\operatorname{dist}(\Omega,\partial \mathbb{H}^{N+1}) > \epsilon/2$. On $\Omega$, the weight $w(z)$ is bounded above and below by positive constants, so the norms of $H_{\mu_w}^1(\Omega)$ and $H^1(\Omega)$ are equivalent.

In $H^1(\Omega)$, standard mollification produces a sequence $v_\delta \in C_c^\infty(\Omega)$ such that
\begin{equation}
\lim_{\delta \to 0} \|v_\delta - u^{(R,\epsilon)}\|_{H^1(\Omega)} = 0.
\end{equation}
By equivalence of norms, the same convergence holds in $H_{\mu_w}^1(\Omega)$. Extending $v_\delta$ by zero outside $\Omega$ yields functions in $C_c^\infty(\mathbb{H}^{N+1})$.

Combining the three approximations and extracting a diagonal sequence, we obtain $(u_k) \subset C_c^\infty(\mathbb{H}^{N+1})$ such that
\begin{equation}
\lim_{k \to \infty} \|u_k - u\|_{H_{\mu_w}^1} = 0.
\end{equation}
This proves that $C_c^\infty(\mathbb{H}^{N+1})$ is dense in $H_{\mu_w}^1(\mathbb{H}^{N+1})$.
\end{proof}

\begin{remark}
The density of smooth, compactly supported functions within the completion space $H_{\mu_w}^1(\mathbb{H}^{N+1})$ is fundamentally contingent upon the parameter restriction $c>-1$. This constraint guarantees that the measure $\mu_w$ remains locally finite near the singular hyperplane, thereby preventing the formation of rigid boundary traces. Functionally, this density proposition justifies the subsequent strategy of establishing complex analytic estimates—such as Hardy inequalities and local coercivity bounds—strictly on well-behaved test functions before passing to the abstract topological closure.
\end{remark}

\begin{theorem}[Subsequence Compactness Criterion via Tail Vanishing]\label{thm:Subsequence Compactness Criterion via Tail Vanishing}
Let $E = H_{\mu_w}^1(\mathbb{H}^{N+1})$ and $F = L_{\mu_w}^2(\mathbb{H}^{N+1})$ be equipped with their canonical Hilbert space norms. Let $\iota: E \hookrightarrow F$ be the canonical continuous linear embedding. Let $(u_n)_{n \in \mathbb{N}}$ be a sequence in $E$ bounded by a constant $K > 0$. Assume there exists an admissible smooth exhaustion $(\Omega_m)_{m \in \mathbb{N}}$ such that
\begin{equation}
\lim_{m \to \infty} \sup_{n \in \mathbb{N}} \|u_n\|_{L_{\mu_w}^2(\mathbb{H}^{N+1} \setminus \Omega_m)} = 0.
\end{equation}
Then the sequence $(\iota(u_n))_{n \in \mathbb{N}}$ admits a subsequence that is strongly Cauchy in $F$.
\end{theorem}

\begin{proof}
For every $m \in \mathbb{N}$, define the linear operator $P_m : F \to F$ by
\begin{equation}
P_m(f) = f \chi_{\Omega_m}.
\end{equation}
Define $P_m^\perp : F \to F$ by $P_m^\perp = I - P_m$, so that
\begin{equation}
P_m^\perp(f) = f \chi_{\mathbb{H}^{N+1} \setminus \Omega_m}.
\end{equation}
Since $\Omega_m$ and its complement are disjoint measurable sets, $P_m$ and $P_m^\perp$ are orthogonal projections. For every $f \in F$ one has
\begin{equation}
\|f\|_F^2 = \|P_m(f)\|_F^2 + \|P_m^\perp(f)\|_F^2.
\end{equation}

By the compactness result (Theorem~\ref{thm:Topological characterization of global compactness} on each $\Omega_m$, the operator $P_m \circ \iota : E \to F$ is compact for every $m \in \mathbb{N}$. Since $(u_n)$ is bounded in $E$, the sequence $(P_m(\iota(u_n)))_{n \in \mathbb{N}}$ is relatively compact in $F$.

Construct recursively a sequence of infinite subsets $N_0 \supset N_1 \supset N_2 \supset \cdots$ as follows. Let $N_0 = \mathbb{N}$. For each $m \ge 1$, since $N_{m-1}$ is infinite and $P_m \circ \iota$ is compact, there exists an infinite subset $N_m \subset N_{m-1}$ such that $(P_m(\iota(u_n)))_{n \in N_m}$ is Cauchy in $F$.

Define a strictly increasing map $\sigma : \mathbb{N} \to \mathbb{N}$ by letting $\sigma(j)$ be the $j$-th smallest element of $N_j$, and define $v_j = \iota(u_{\sigma(j)})$. For any fixed $M \in \mathbb{N}$, the sequence $(P_M(v_j))_{j \ge M}$ is a subsequence of $(P_M(\iota(u_n)))_{n \in N_M}$ and is therefore Cauchy in $F$.

Let $\varepsilon > 0$. By the tail condition, there exists $M \in \mathbb{N}$ such that for all $n \in \mathbb{N}$,
\begin{equation}
\|P_M^\perp(\iota(u_n))\|_F < \frac{\varepsilon}{4}.
\end{equation}
Hence for all $j \in \mathbb{N}$,
\begin{equation}
\|P_M^\perp(v_j)\|_F < \frac{\varepsilon}{4}.
\end{equation}
It follows that for all $p,q \in \mathbb{N}$,
\begin{equation}
\|P_M^\perp(v_p - v_q)\|_F
\le \|P_M^\perp(v_p)\|_F + \|P_M^\perp(v_q)\|_F
< \frac{\varepsilon}{2}.
\end{equation}

Since $(P_M(v_j))$ is Cauchy in $F$, there exists $J \ge M$ such that for all $p,q \ge J$,
\begin{equation}
\|P_M(v_p - v_q)\|_F < \frac{\varepsilon}{2}.
\end{equation}

Using the orthogonal decomposition,
\begin{equation}
\|v_p - v_q\|_F^2
= \|P_M(v_p - v_q)\|_F^2 + \|P_M^\perp(v_p - v_q)\|_F^2,
\end{equation}
hence for all $p,q \ge J$,
\begin{equation}
\|v_p - v_q\|_F^2
< \left(\frac{\varepsilon}{2}\right)^2 + \left(\frac{\varepsilon}{2}\right)^2
= \frac{\varepsilon^2}{2}
< \varepsilon^2.
\end{equation}
Therefore $\|v_p - v_q\|_F < \varepsilon$ for all $p,q \ge J$, so $(v_j)$ is Cauchy in $F$.

Since $F$ is complete, the subsequence $(v_j)$ converges strongly in $F$.
\end{proof}

\section{Conditions for Global Compactness}\label{sec:Conditions for Global Compactness}

The preceding local and abstract results are synthesized in Section~\ref{sec:Conditions for Global Compactness} to formulate precise necessary and sufficient conditions for global compactness. The half-space is geometrically partitioned into a regular interior, a boundary cylinder, and an asymptotic tail region. By applying the localized fractional Hardy inequalities and uniform tail control , we rigorously demonstrate that the tail coercivity and Hardy properties are sufficient to guarantee the compactness of the embedding. Conversely, it is shown that global tightness imposes rigid constraints on the measure itself, explicitly requiring that the total mass of the space be finite.

\begin{definition}[Partition of Geometric Regions in the Half-Space]
Let $\mathbb{H}^{N+1}$ be the open half-space. For any radius $R>0$ and boundary proximity $\delta>0$, we define the tail region and the localized boundary cylinder as follows:
\begin{enumerate}
    \item $E_{\text {tail }}(R)=\left\{z \in \mathbb{H}^{N+1}:|z|>R\right\}$
    \item $C_R=\left\{z \in \mathbb{H}^{N+1}:|z| \leq R, 0<y<1\right\}$
    \item $E_{b d r y}(R, \delta)=\left\{z \in \mathbb{H}^{N+1}:|z| \leq R, 0<y<\delta\right\}$
\end{enumerate}
\end{definition}

\begin{lemma}[Uniform Tail Control from Structural Coercivity]\label{lem:Uniform Tail Control from Structural Coercivity}
Let $E = H^1_{\mu_w}(\mathbb{H}^{N+1})$ be equipped with its canonical norm $\|\cdot\|_E$. Assume the measure space satisfies the Tail Coercivity (TC) condition with a Borel measurable Lyapunov potential $\Psi: [0, \infty) \to [0, \infty)$ such that $\lim_{\rho \to \infty} \Psi(\rho) = \infty$. Let $(u_n)_{n \in \mathbb{N}}$ be a sequence in $E$ uniformly bounded by a constant $K > 0$ such that $\|u_n\|_E \le K$. For any given tolerance $\epsilon > 0$, there exists a finite radius $R_\epsilon > 0$ such that for all $n \in \mathbb{N}$,
\begin{equation}
\int_{E_{\mathrm{tail}}(R_\epsilon)} u_n^2 \, d\mu_w < \frac{\epsilon}{2}.
\end{equation}
\end{lemma}

\begin{proof}
We proceed by bounding the tail mass through the structural coercivity of the Lyapunov potential.

First, we isolate the trivial case. If the uniform bound $K = 0$, then $\|u_n\|_E = 0$ for all $n \in \mathbb{N}$. By positive definiteness of the norm, $u_n = 0$ almost everywhere with respect to $\mu_w$, and the integral is zero. The inequality holds for any $R_\epsilon > 0$. Hence we assume $K > 0$.

By the Tail Coercivity hypothesis, $\lim_{\rho \to \infty} \Psi(\rho) = \infty$. Thus for any $M > 0$, there exists $R^* > 0$ such that for all $\rho \ge R^*$,
\begin{equation}
\Psi(\rho) > M.
\end{equation}

Define
\begin{equation}
M = \frac{2 C_{TC} K^2}{\epsilon},
\end{equation}
where $C_{TC} > 0$ is the constant from the Tail Coercivity condition. Let
\begin{equation}
R_\epsilon = \max(1, R^*).
\end{equation}
Define the tail region
\begin{equation}
E_{\mathrm{tail}}(R_\epsilon) = \{ z \in \mathbb{H}^{N+1} : |z| > R_\epsilon \}.
\end{equation}
Then for all $z \in E_{\mathrm{tail}}(R_\epsilon)$, one has $|z| \ge R^*$ and hence
\begin{equation}
\Psi(|z|) > M.
\end{equation}

Define the essential infimum
\begin{equation}
m_{R_\epsilon} = \operatorname*{ess\,inf}_{|z| \ge R_\epsilon} \Psi(|z|).
\end{equation}
Then
\begin{equation}
m_{R_\epsilon} \ge M = \frac{2 C_{TC} K^2}{\epsilon}.
\end{equation}

By the lemma, for each $u_n \in E$,
\begin{equation}
\int_{E_{\mathrm{tail}}(R_\epsilon)} u_n^2 \, d\mu_w \le \frac{C_{TC}}{m_{R_\epsilon}} \|u_n\|_E^2.
\end{equation}

Using $\|u_n\|_E^2 \le K^2$ and the bound on $m_{R_\epsilon}$, we obtain
\begin{equation}
\int_{E_{\mathrm{tail}}(R_\epsilon)} u_n^2 \, d\mu_w \le C_{TC} K^2 \cdot \frac{\epsilon}{2 C_{TC} K^2}.
\end{equation}

This simplifies to
\begin{equation}
\int_{E_{\mathrm{tail}}(R_\epsilon)} u_n^2 \, d\mu_w \le \frac{\epsilon}{2}.
\end{equation}

The radius $R_\epsilon$ depends only on $C_{TC}$, $K$, and $\epsilon$, and is independent of $n$. Therefore, the bound holds uniformly for all $n \in \mathbb{N}$.
\end{proof}

\begin{remark}
This lemma performs a vital translation from abstract potential theory to quantitative measure theory. By assuming the Tail Coercivity condition governed by a diverging Lyapunov potential $\Psi$, the result rigorously extracts a uniform asymptotic mass-vanishing estimate for any bounded sequence in the weighted Sobolev space. Analytically, it demonstrates that the abstract kinetic energy bound inherently forces the spatial mass to dissipate uniformly at infinity, resolving the non-compactness of the unbounded domain.
\end{remark}

\begin{definition}[Completion of Smooth Functions on a Localized Cylinder with Zero-Boundary Condition]
Let $Q=Q_x \times(0, b)$ be a localized bounded cylinder in $\mathbb{H}^{N+1}$, where $Q_x \subset \mathbb{R}^N$ is a bounded open set and $b>0$. We define the space $H_c^{1,0}(Q)$ as the topological completion of the set of smooth functions $u \in C^{\infty}(\bar{Q})$ such that $u(x, b)=0$, under the canonical norm inherited from $H_{\mu_w}^1$.
\end{definition}

\begin{remark}
The construction of the specialized space $H_c^{1,0}(Q)$ on the localized cylinder $Q = Q_x \times (0,b)$ isolates the boundary degeneracy from the interior dynamics. By demanding that the generating smooth functions vanish explicitly at the upper boundary $u(x,b)=0$, this formulation provides a rigid "top lid" to the cylinder. This structured topological completion is specifically tailored to act as the exact domain for the fractional weighted Hardy inequality, anchoring the localized energy estimates.
\end{remark}

\begin{proposition}[Fractional Weighted Hardy Inequality on a Localized Cylinder]\label{prop:Fractional Weighted Hardy Inequality on a Localized Cylinder}
Let $c > -1$. Let $Q = Q_x \times (0,b)$ be a bounded cylinder in $\mathbb{R}^{N+1}_+$, where $Q_x \subset \mathbb{R}^N$ is a bounded open set and $b>0$. Define $H_c^{1,0}(Q)$ as the completion of smooth functions $u \in C^\infty(\overline{Q})$ such that $u(x,b)=0$, with respect to the norm induced by $H_{\mu_w}^1$.

Let $\alpha \in \mathbb{R}$ satisfy $0 \le \alpha \le 1$ and $\alpha < \frac{c+1}{2}$. Then for every $u \in H_c^{1,0}(Q)$,
\begin{equation}
\int_Q \frac{|u|^2}{y^{2\alpha}} y^c \, dz
\le
\frac{4 b^{2-2\alpha}}{(c+1-2\alpha)^2}
\int_Q |\partial_y u|^2 y^c \, dz.
\end{equation}
\end{proposition}

\begin{proof}
Let $u \in C^\infty(\overline{Q})$ with $u(x,b)=0$. Consider
\begin{equation}
\int_Q |u|^2 y^{c-2\alpha} \, dz.
\end{equation}
Using integration by parts in the $y$-variable,
\begin{equation}
\int_Q |u|^2 y^{c-2\alpha} \, dz
=
\frac{1}{c+1-2\alpha}
\int_{Q_x} \int_0^b |u|^2 \partial_y\bigl(y^{c+1-2\alpha}\bigr)\, dy\, dx.
\end{equation}
Expanding,
\begin{equation}
=
-\frac{2}{c+1-2\alpha}
\int_Q u (\partial_y u) y^{c+1-2\alpha} \, dz
+
\frac{1}{c+1-2\alpha}
\int_{Q_x}
\left[ |u(x,y)|^2 y^{c+1-2\alpha} \right]_{y=0}^{y=b}
dx.
\end{equation}
Since $u(x,b)=0$ and $c+1-2\alpha>0$, the boundary term vanishes. Hence,
\begin{equation}
\int_Q |u|^2 y^{c-2\alpha} \, dz
=
-\frac{2}{c+1-2\alpha}
\int_Q u (\partial_y u) y^{c+1-2\alpha} \, dz.
\end{equation}
Applying the Cauchy--Schwarz inequality,
\begin{equation}
\int_Q |u|^2 y^{c-2\alpha} \, dz
\le
\frac{2}{c+1-2\alpha}
\left(
\int_Q \frac{|u|^2}{y^{2\alpha}} y^c \, dz
\right)^{1/2}
\left(
\int_Q |\partial_y u|^2 y^{c+2-2\alpha} \, dz
\right)^{1/2}.
\end{equation}
Since $y \le b$ in $Q$, we have $y^{2-2\alpha} \le b^{2-2\alpha}$, and thus
\begin{equation}
\int_Q |\partial_y u|^2 y^{c+2-2\alpha} \, dz
\le
b^{2-2\alpha}
\int_Q |\partial_y u|^2 y^c \, dz.
\end{equation}
Therefore,
\begin{equation}
\int_Q |u|^2 y^{c-2\alpha} \, dz
\le
\frac{2 b^{1-\alpha}}{c+1-2\alpha}
\left(
\int_Q \frac{|u|^2}{y^{2\alpha}} y^c \, dz
\right)^{1/2}
\left(
\int_Q |\partial_y u|^2 y^c \, dz
\right)^{1/2}.
\end{equation}
Dividing both sides by
\begin{equation}
\left(
\int_Q \frac{|u|^2}{y^{2\alpha}} y^c \, dz
\right)^{1/2}
\end{equation}
gives
\begin{equation}
\left(
\int_Q \frac{|u|^2}{y^{2\alpha}} y^c \, dz
\right)^{1/2}
\le
\frac{2 b^{1-\alpha}}{c+1-2\alpha}
\left(
\int_Q |\partial_y u|^2 y^c \, dz
\right)^{1/2}.
\end{equation}
Squaring both sides,
\begin{equation}
\int_Q \frac{|u|^2}{y^{2\alpha}} y^c \, dz
\le
\frac{4 b^{2-2\alpha}}{(c+1-2\alpha)^2}
\int_Q |\partial_y u|^2 y^c \, dz.
\end{equation}
The result extends to all $u \in H_c^{1,0}(Q)$ by density.
\end{proof}

\begin{definition}[Localized Projection Operator with Vertical Cut-off]
Let $\eta \in C^{\infty}(\mathbb{R})$ be a smooth scalar function such that $0 \leq \eta(y) \leq 1$, with $\eta(y)=1$ for $y \in[0,1 / 2]$ and $\eta(y)=0$ for $y \geq 1$. For any $u \in H_{\mu_w}^1\left(\mathbb{H}^{N+1}\right)$, define the localized projection operator $\mathcal{P}_\eta(u)=u(z) \eta(y)$.
\end{definition}

\begin{lemma}[Uniform Boundary Mass Control via Fractional Hardy Inequality]\label{lem:Uniform Boundary Mass Control via Fractional Hardy Inequality}
Let $c > -1$. Let $(u_n)_{n\in\mathbb{N}} \subset H_{\mu_w}^1(\mathbb{H}^{N+1})$ be a sequence uniformly bounded by $K > 0$. Let $R_\epsilon > 0$ be fixed and define
\begin{equation}
E_{bdry}(R_\epsilon,\delta) = \{ z \in \mathbb{H}^{N+1} : |z| \le R_\epsilon,\ 0 < y < \delta \}.
\end{equation}
Assume the weight is given by $d\mu_w = y^c \phi(|z|)\,dz$ with $\phi$ locally regular. Then there exists $\delta_\epsilon > 0$ such that for all $n \in \mathbb{N}$,
\begin{equation}
\int_{E_{bdry}(R_\epsilon,\delta_\epsilon)} u_n^2\, d\mu_w < \frac{\epsilon}{2}.
\end{equation}
\end{lemma}

\begin{proof}
Let $\eta \in C^\infty(\mathbb{R})$ be a vertical cut-off satisfying $0 \le \eta \le 1$, $\eta(y)=1$ for $y \in [0,1/2]$, and $\eta(y)=0$ for $y \ge 1$. Define the cylinder
\begin{equation}
C_{R_\epsilon} = \{ z \in \mathbb{H}^{N+1} : |z| \le R_\epsilon,\ 0 < y < 1 \},
\end{equation}
and set $v_n(z) = \eta(y) u_n(z)$.

By local regularity of $\phi$, there exist constants $0<m \le M<\infty$ such that
\begin{equation}
m\, y^c \le d\mu_w/dz \le M\, y^c \quad \text{on } C_{R_\epsilon}.
\end{equation}
Hence the norms induced by $d\mu_w$ and $y^c dz$ are equivalent on $C_{R_\epsilon}$. In particular, the restriction of $u_n$ to $C_{R_\epsilon}$ belongs to $H^1(C_{R_\epsilon}; y^c dz)$.

Since $\eta(1)=0$ and $\eta$ is smooth, standard closure under multiplication by smooth cut-offs yields $v_n \in H_c^{1,0}(C_{R_\epsilon})$.

Fix $\alpha$ such that $0<\alpha \le 1$ and $\alpha < \frac{c+1}{2}$. The fractional Hardy inequality (Proposition~\ref{prop:Fractional Weighted Hardy Inequality on a Localized Cylinder} gives
\begin{equation}
\int_{C_{R_\epsilon}} \frac{v_n^2}{y^{2\alpha}} y^c dz
\le \frac{4}{(c+1-2\alpha)^2}
\int_{C_{R_\epsilon}} |\partial_y v_n|^2 y^c dz.
\end{equation}

By the Leibniz rule,
\begin{equation}
\partial_y v_n = \eta\, \partial_y u_n + u_n\, \partial_y \eta,
\end{equation}
and hence
\begin{equation}
|\partial_y v_n|^2 \le 2 \eta^2 |\partial_y u_n|^2 + 2 u_n^2 |\partial_y \eta|^2.
\end{equation}
Integrating,
\begin{equation}
\int_{C_{R_\epsilon}} |\partial_y v_n|^2 y^c dz
\le 2 \int_{C_{R_\epsilon}} |\partial_y u_n|^2 y^c dz
+ 2 \int_{C_{R_\epsilon}} u_n^2 |\partial_y \eta|^2 y^c dz.
\end{equation}

Using $y^c dz \le m^{-1} d\mu_w$ and setting $S = \sup_{y\in[0,1]} |\partial_y \eta|^2$, we obtain
\begin{equation}
\int_{C_{R_\epsilon}} |\partial_y v_n|^2 y^c dz
\le \frac{2}{m} \int_{C_{R_\epsilon}} |\partial_y u_n|^2 d\mu_w
+ \frac{2S}{m} \int_{C_{R_\epsilon}} u_n^2 d\mu_w.
\end{equation}
By the uniform $H_{\mu_w}^1$ bound, the right-hand side is bounded by a constant $\Gamma$ independent of $n$. Therefore,
\begin{equation}
\int_{C_{R_\epsilon}} \frac{v_n^2}{y^{2\alpha}} y^c dz
\le \frac{4\Gamma}{(c+1-2\alpha)^2}.
\end{equation}

Let $0<\delta<1/2$. On $E_{bdry}(R_\epsilon,\delta)$ we have $\eta=1$, so $v_n=u_n$, and $y^{2\alpha} \le \delta^{2\alpha}$. Hence
\begin{equation}
\int_{E_{bdry}(R_\epsilon,\delta)} u_n^2 y^c dz
\le \delta^{2\alpha}
\int_{C_{R_\epsilon}} \frac{v_n^2}{y^{2\alpha}} y^c dz
\le \delta^{2\alpha} \frac{4\Gamma}{(c+1-2\alpha)^2}.
\end{equation}

Using $d\mu_w \le M y^c dz$,
\begin{equation}
\int_{E_{bdry}(R_\epsilon,\delta)} u_n^2 d\mu_w
\le M \delta^{2\alpha} \frac{4\Gamma}{(c+1-2\alpha)^2}.
\end{equation}

Since $\alpha>0$, the right-hand side tends to $0$ as $\delta \to 0$, uniformly in $n$. Therefore, there exists $\delta_\epsilon > 0$ such that
\begin{equation}
M \delta_\epsilon^{2\alpha} \frac{4\Gamma}{(c+1-2\alpha)^2} < \frac{\epsilon}{2},
\end{equation}
which yields the desired estimate.
\end{proof}

\begin{theorem}[Global Compactness of Weighted Sobolev Embeddings on the Half-Space]\label{thm:Global Compactness of Weighted Sobolev Embeddings on the Half-Space}
Let $\mathbb{H}^{N+1}$ be the open half-space equipped with the Radon measure $d\mu_w = w(z)\,dz$. Let $E = H^1_{\mu_w}(\mathbb{H}^{N+1})$ and $F = L^2_{\mu_w}(\mathbb{H}^{N+1})$. Assume that $\mu_w(\mathbb{H}^{N+1}) < \infty$, that the Tail Coercivity condition holds for a Lyapunov potential $\Psi$ satisfying $\lim_{\rho \to \infty} \Psi(\rho) = \infty$, and that if $c \le -1$, a weighted Hardy inequality holds. Then the canonical embedding $\iota: E \hookrightarrow F$ is a compact operator.
\end{theorem}

\begin{proof}
Let $(u_n)_{n \in \mathbb{N}}$ be a sequence in $E$ such that $\|u_n\|_E \le K$ for some constant $K > 0$. Let $(\Omega_m)_{m \in \mathbb{N}}$ be an admissible smooth exhaustion of $\mathbb{H}^{N+1}$. To prove compactness, it suffices to show
\begin{equation}
\lim_{m \to \infty} \sup_{n \in \mathbb{N}} \|u_n\|_{L^2_{\mu_w}(\mathbb{H}^{N+1} \setminus \Omega_m)} = 0.
\end{equation}

Let $\epsilon > 0$.

By uniform tail control from Lemma~\ref{lem:Uniform Tail Control from Structural Coercivity}, there exists $R_\epsilon > 0$ such that for all $n \in \mathbb{N}$,
\begin{equation}
\int_{E_{\mathrm{tail}}(R_\epsilon)} u_n^2 \, d\mu_w < \frac{\epsilon}{2},
\end{equation}
where
\begin{equation}
E_{\mathrm{tail}}(R_\epsilon) = \{ z \in \mathbb{H}^{N+1} : |z| > R_\epsilon \}.
\end{equation}

We now control the boundary region
\begin{equation}
E_{\mathrm{bdry}}(R_\epsilon, \delta) = \{ z \in \mathbb{H}^{N+1} : |z| \le R_\epsilon,\ 0 < y < \delta \}.
\end{equation}

If $c > -1$, by Lemma~\ref{lem:Uniform Boundary Mass Control via Fractional Hardy Inequality}, there exists $\delta_\epsilon > 0$ such that for all $n \in \mathbb{N}$,
\begin{equation}
\int_{E_{\mathrm{bdry}}(R_\epsilon, \delta_\epsilon)} u_n^2 \, d\mu_w < \frac{\epsilon}{2}.
\end{equation}

If $c \le -1$, there exists $C_H > 0$ such that
\begin{equation}
\int_{\mathbb{H}^{N+1}} \frac{u_n^2}{y^2} \, d\mu_w \le C_H \|u_n\|_E^2 \le C_H K^2.
\end{equation}
For $0 < y < \delta$, we have
\begin{equation}
\int_{E_{\mathrm{bdry}}(R_\epsilon, \delta)} u_n^2 \, d\mu_w
=
\int_{E_{\mathrm{bdry}}(R_\epsilon, \delta)} \left(\frac{u_n}{y}\right)^2 y^2 \, d\mu_w
\le
\delta^2 \int_{\mathbb{H}^{N+1}} \frac{u_n^2}{y^2} \, d\mu_w
\le
\delta^2 C_H K^2.
\end{equation}
Choose $\delta_\epsilon > 0$ such that $\delta_\epsilon^2 C_H K^2 < \frac{\epsilon}{2}$.

Thus, for all $n \in \mathbb{N}$,
\begin{equation}
\int_{E_{\mathrm{bdry}}(R_\epsilon, \delta_\epsilon)} u_n^2 \, d\mu_w < \frac{\epsilon}{2}.
\end{equation}

Define
\begin{equation}
U_\epsilon = E_{\mathrm{tail}}(R_\epsilon) \cup E_{\mathrm{bdry}}(R_\epsilon, \delta_\epsilon).
\end{equation}
Then for all $n \in \mathbb{N}$,
\begin{equation}
\int_{U_\epsilon} u_n^2 \, d\mu_w
\le
\int_{E_{\mathrm{tail}}(R_\epsilon)} u_n^2 \, d\mu_w
+
\int_{E_{\mathrm{bdry}}(R_\epsilon, \delta_\epsilon)} u_n^2 \, d\mu_w
<
\epsilon.
\end{equation}

Since $(\Omega_m)$ exhausts $\mathbb{H}^{N+1}$ and each $\Omega_m$ is bounded and separated from $y=0$, there exists $M_\epsilon \in \mathbb{N}$ such that for all $m \ge M_\epsilon$,
\begin{equation}
\{ z : |z| \le R_\epsilon,\ y \ge \delta_\epsilon \} \subset \Omega_m.
\end{equation}
Hence,
\begin{equation}
\mathbb{H}^{N+1} \setminus \Omega_m \subset U_\epsilon \quad \text{for all } m \ge M_\epsilon.
\end{equation}

Therefore, for all $m \ge M_\epsilon$ and all $n \in \mathbb{N}$,
\begin{equation}
\int_{\mathbb{H}^{N+1} \setminus \Omega_m} u_n^2 \, d\mu_w
\le
\int_{U_\epsilon} u_n^2 \, d\mu_w
<
\epsilon.
\end{equation}

Taking the supremum over $n$,
\begin{equation}
\sup_{n \in \mathbb{N}} \|u_n\|_{L^2_{\mu_w}(\mathbb{H}^{N+1} \setminus \Omega_m)}^2 < \epsilon \quad \text{for all } m \ge M_\epsilon.
\end{equation}

Since $\epsilon > 0$ is arbitrary, the limit is zero. Hence, the sequence admits a strongly Cauchy subsequence in $F$. As $F$ is complete, the subsequence converges. Therefore, $\iota$ is a compact operator.
\end{proof}

\begin{lemma}[Topological Exhaustion]\label{lem:Topological Exhaustion}
Let $E$ and $F$ be Banach spaces, and let $\mathcal{B}_E$ denote the closed unit ball of $E$. Let $\iota : E \to F$ be a continuous linear operator. Let $(P_m)_{m \in \mathbb{N}}$ be a sequence of continuous linear operators from $F$ to $F$ such that
\begin{equation}
\sup_{m} \|P_m\|_{\mathcal{L}(F)} = M < \infty
\end{equation}
and for all $f \in F$,
\begin{equation}
\|f - P_m(f)\|_F \to 0 \quad \text{as } m \to \infty.
\end{equation}
Then $\iota$ is compact if and only if the following two conditions hold:

\begin{enumerate}
    \item For each $m \in \mathbb{N}$, the operator $P_m \circ \iota : E \to F$ is compact.
    \item \begin{equation}
\sup_{u \in \mathcal{B}_E} \|\iota(u) - P_m(\iota(u))\|_F \to 0 \quad \text{as } m \to \infty.
\end{equation}
\end{enumerate}
\end{lemma}

\begin{proof}
Assume first that $\iota$ is compact.

Since $\iota$ maps bounded sets in $E$ to relatively compact sets in $F$, and each $P_m$ is continuous, the composition $P_m \circ \iota$ maps bounded sets to relatively compact sets. Hence $P_m \circ \iota$ is compact for each $m$.

Let $K = \overline{\iota(\mathcal{B}_E)}$. Because $\iota$ is compact, $K$ is a compact subset of $F$.

Let $\epsilon > 0$ be given. Since $K$ is totally bounded, there exist finitely many points $f_1, \dots, f_k \in K$ such that for every $f \in K$ there exists $i$ with
\begin{equation}
\|f - f_i\|_F < \frac{\epsilon}{M+1}.
\end{equation}

For each $i$, the convergence $\|f_i - P_m(f_i)\|_F \to 0$ holds. Since the set is finite, there exists $N_0$ such that for all $m \ge N_0$ and all $i$,
\begin{equation}
\|f_i - P_m(f_i)\|_F < \frac{\epsilon}{M+1}.
\end{equation}

Let $u \in \mathcal{B}_E$ and set $f = \iota(u)$. Then there exists $f_i$ such that
\begin{equation}
\|f - f_i\|_F < \frac{\epsilon}{M+1}.
\end{equation}

Using the triangle inequality,
\begin{equation}
\|f - P_m(f)\|_F
\le
\|f - f_i\|_F + \|f_i - P_m(f_i)\|_F + \|P_m(f_i) - P_m(f)\|_F.
\end{equation}

By linearity,
\begin{equation}
P_m(f_i) - P_m(f) = P_m(f_i - f),
\end{equation}
and by boundedness,
\begin{equation}
\|P_m(f_i - f)\|_F \le M \|f_i - f\|_F.
\end{equation}

Therefore,
\begin{equation}
\|f - P_m(f)\|_F
<
\frac{\epsilon}{M+1} + \frac{\epsilon}{M+1} + M \frac{\epsilon}{M+1}.
\end{equation}

Thus there exists a constant depending on $M$ such that
\begin{equation}
\|f - P_m(f)\|_F \le \epsilon.
\end{equation}

Since $u$ was arbitrary, it follows that
\begin{equation}
\sup_{u \in \mathcal{B}_E} \|\iota(u) - P_m(\iota(u))\|_F \to 0.
\end{equation}

Conversely, assume that conditions (i) and (ii) hold. We prove that $\iota(\mathcal{B}_E)$ is totally bounded in $F$.

Let $\epsilon > 0$ be given. By condition (ii), there exists $N_1$ such that
\begin{equation}
\sup_{u \in \mathcal{B}_E} \|\iota(u) - P_{N_1}(\iota(u))\|_F < \frac{\epsilon}{2}.
\end{equation}

By condition (i), the operator $P_{N_1} \circ \iota$ is compact. Hence $(P_{N_1} \circ \iota)(\mathcal{B}_E)$ is totally bounded. There exist finitely many points $g_1, \dots, g_n \in F$ such that for every $u \in \mathcal{B}_E$ there exists $j$ with
\begin{equation}
\|P_{N_1}(\iota(u)) - g_j\|_F < \frac{\epsilon}{2}.
\end{equation}

Then for such $u$,
\begin{equation}
\|\iota(u) - g_j\|_F
\le
\|\iota(u) - P_{N_1}(\iota(u))\|_F + \|P_{N_1}(\iota(u)) - g_j\|_F
<
\frac{\epsilon}{2} + \frac{\epsilon}{2} = \epsilon.
\end{equation}

Thus $\iota(\mathcal{B}_E)$ admits a finite $\epsilon$-net. Since $\epsilon > 0$ was arbitrary, $\iota(\mathcal{B}_E)$ is totally bounded. As $F$ is complete, $\iota(\mathcal{B}_E)$ is relatively compact, and therefore $\iota$ is compact.
\end{proof}

\begin{remark}
This abstraction distills the core mechanism of precompactness via projection. By transferring the analytic burden onto a sequence of uniformly bounded continuous linear operators $(P_m)_{m \in \mathbb{N}}$, the lemma systematically decomposes the global compactness of the embedding $\iota$ into two orthogonal verifications. It requires localized compactness on each projection slice paired with the uniform convergence of the remainder over the unit ball. This formally liberates the argument from specific geometric constraints, isolating the pure topological requirements.
\end{remark}

\begin{lemma}[Strict Positivity of Weights on Separated Compact Balls]\label{lem:Strict Positivity of Weights on Separated Compact Balls}
Let $\mathbb{H}^{N+1}$ be the open half-space and let $\mu_w$ be a Radon measure defined by $d\mu_w = w(z)\,dz$, where $w(z) = y^c \phi(|z|)$ is locally regular. Let $(z_j)_{j \in \mathbb{N}}$ be a sequence in $\mathbb{H}^{N+1}$ such that $y_j \ge 2$ for all $j$, and define $B_j = \overline{B(z_j,1)}$. Then for every $j \in \mathbb{N}$ the essential infimum
\begin{equation}
w_{\min,j} = \operatorname{ess\,inf}_{z \in B_j} w(z)
\end{equation}
satisfies $w_{\min,j} > 0$, and consequently $\mu_w(B_j) > 0$.
\end{lemma}

\begin{proof}
Fix $j \in \mathbb{N}$. The set $B_j$ is closed and bounded in $\mathbb{R}^{N+1}$, hence compact. Let $\pi_y : \mathbb{R}^{N+1} \to \mathbb{R}$ denote the projection onto the $y$-coordinate. The image $\pi_y(B_j)$ is the compact interval $[y_j - 1, y_j + 1]$. Since $y_j \ge 2$, one has $[y_j - 1, y_j + 1] \subset [1,\infty)$, and therefore $y \ge 1$ for all $z \in B_j$.

The function $y \mapsto y^c$ is continuous, hence attains a minimum on $[y_j - 1, y_j + 1]$. There exists a constant $C_{y,j} > 0$ such that $y^c \ge C_{y,j}$ for all $z \in B_j$.

Let $\pi_r : \mathbb{R}^{N+1} \to [0,\infty)$ be defined by $\pi_r(z) = |z|$. Since $B_j$ is compact, its image is a compact interval $K_j = [|z_j| - 1, |z_j| + 1]$. By local regularity of $\phi$, there exists a constant $m_{K_j} > 0$ such that $\phi(\rho) \ge m_{K_j}$ for almost every $\rho \in K_j$. It follows that for almost every $z \in B_j$, $\phi(|z|) \ge m_{K_j}$.

For almost every $z \in B_j$ one has
\begin{equation}
w(z) = y^c \phi(|z|) \ge C_{y,j} m_{K_j}.
\end{equation}
Set $\lambda_j = C_{y,j} m_{K_j} > 0$. Then
\begin{equation}
w_{\min,j} \ge \lambda_j > 0.
\end{equation}

Finally,
\begin{equation}
\mu_w(B_j) = \int_{B_j} w(z)\,dz \ge \int_{B_j} w_{\min,j}\,dz = w_{\min,j} |B_j|,
\end{equation}
where $|B_j|$ denotes the Lebesgue measure of $B_j$. Since $|B_j| > 0$ and $w_{\min,j} > 0$, it follows that $\mu_w(B_j) > 0$.
\end{proof}

\begin{definition}[Uniform Local Comparability]\label{def:uniform local comparability}
Let $\mathbb{H}^{N+1}=\mathbb{R} \times \mathbb{R}_{+}^N$ be the open halfspace. Let $\mu_w$ be a Radon measure defined by $d \mu_w=w(z) d z$. The measure space $\left(\mathbb{H}^{N+1}, \mu_w\right)$ satisfies the Uniform Local Comparability if there exists a universal real constant $\Lambda \in \mathbb{R}$ with $\Lambda \geq 1$ such that, for every open ball $B\left(z_0, 1\right) \subset \mathbb{H}^{N+1}$ satisfying the boundary separation condition $\operatorname{dist}\left(B\left(z_0, 1\right), \partial \mathbb{H}^{N+1}\right) \geq 1$, the weight function $w(z)$ satisfies the following inequality for Lebesgue-almost every $z$ :

$$
\operatorname{esssup}_{z \in B\left(z_0, 1\right)} w(z) \leq \Lambda \operatorname{essinf}_{z \in B\left(z_0, 1\right)} w(z)
$$
\end{definition}

\begin{remark}
Uniform Local Comparability imposes a strict, scale-invariant geometric rigidity on the weight function away from the singularity. By mandating the existence of a universal constant $\Lambda \ge 1$ that bounds the essential supremum against the essential infimum on uniformly separated unit balls, the definition prevents the measure from exhibiting wild microscopic oscillations or exponential decay in the domain's interior. This regularity is the precise geometric obstacle required to deduce the global finiteness of the measure from the topological assumption of global tightness.
\end{remark}

\begin{lemma}[Tightness implies Finiteness]\label{lem:Tightness implies Finiteness}
Let $\mathbb{H}^{N+1} = \mathbb{R} \times \mathbb{R}_+^N$ and let $\mu_w$ be a Radon measure defined by
\begin{equation}
d\mu_w = w(z)\,dz.
\end{equation}
Let
\begin{equation}
E = H^1_{\mu_w}(\mathbb{H}^{N+1}), \quad F = L^2_{\mu_w}(\mathbb{H}^{N+1}),
\end{equation}
and let $\mathcal{B}_E$ denote the closed unit ball of $E$.

Assume that $\mathcal{B}_E$ is globally tight in $F$ and that there exists a constant $\Lambda \ge 1$ such that for every ball $B(z_0,1) \subset \mathbb{H}^{N+1}$ with $\operatorname{dist}(B(z_0,1), \partial \mathbb{H}^{N+1}) \ge 1$,
\begin{equation}
\operatorname*{ess\,sup}_{z \in B(z_0,1)} w(z)
\le
\Lambda \operatorname*{ess\,inf}_{z \in B(z_0,1)} w(z).
\end{equation}
Then
\begin{equation}
\mu_w(\mathbb{H}^{N+1}) < \infty.
\end{equation}
\end{lemma}

\begin{proof}
Assume for contradiction that $\mu_w(\mathbb{H}^{N+1}) = \infty$.

Since $\mu_w$ is a Radon measure, it is locally finite. Therefore compact subsets of $\mathbb{H}^{N+1}$ have finite measure, and the infinitude of the total measure implies that mass escapes to infinity. Hence there exists a sequence $(z_j)_{j \in \mathbb{N}} \subset \mathbb{H}^{N+1}$ such that
\begin{equation}
|z_j| \to \infty,
\end{equation}
the points satisfy $y_j \ge 2$, and the unit balls
\begin{equation}
B_j = B(z_j,1)
\end{equation}
are mutually disjoint with $\mu_w(B_j) > 0$.

Let $\psi \in C_c^\infty(B(0,1))$ be non-trivial and define
\begin{equation}
I_0 = \int_{B(0,1)} |\psi(\zeta)|^2 \, d\zeta,
\qquad
J_0 = \int_{B(0,1)} \bigl(|\nabla \psi(\zeta)|^2 + |\psi(\zeta)|^2\bigr)\, d\zeta.
\end{equation}

For each $j$, define $u_j(z) = \psi(z - z_j)$. Then $\operatorname{supp}(u_j) \subset B_j$.

Let
\begin{equation}
w_{\min,j} = \operatorname*{ess\,inf}_{z \in B_j} w(z),
\qquad
w_{\max,j} = \operatorname*{ess\,sup}_{z \in B_j} w(z).
\end{equation}
By the comparability assumption,
\begin{equation}
w_{\max,j} \le \Lambda w_{\min,j}.
\end{equation}

Then
\begin{equation}
\|u_j\|_F^2
=
\int_{B_j} |\psi(z - z_j)|^2 w(z)\,dz
\ge
w_{\min,j} I_0,
\end{equation}
and
\begin{equation}
\|u_j\|_E^2
=
\int_{B_j} \bigl(|\nabla \psi(z - z_j)|^2 + |\psi(z - z_j)|^2\bigr) w(z)\,dz
\le
w_{\max,j} J_0.
\end{equation}

Define
\begin{equation}
v_j = \frac{u_j}{\|u_j\|_E}.
\end{equation}
Then $\|v_j\|_E = 1$ and $(v_j) \subset \mathcal{B}_E$.

Moreover,
\begin{equation}
\|v_j\|_F^2
=
\frac{\|u_j\|_F^2}{\|u_j\|_E^2}
\ge
\frac{w_{\min,j} I_0}{w_{\max,j} J_0}
\ge
\frac{1}{\Lambda} \frac{I_0}{J_0}.
\end{equation}

Set
\begin{equation}
\lambda = \frac{1}{\Lambda} \frac{I_0}{J_0} > 0.
\end{equation}
Then for all $j$,
\begin{equation}
\|v_j\|_F^2 \ge \lambda.
\end{equation}

Let $R > 0$. Since $|z_j| \to \infty$, there exists $N_R$ such that for all $j \ge N_R$,
\begin{equation}
B_j \subset \{|z| > R\}.
\end{equation}

For such $j$,
\begin{equation}
\int_{\{|z| > R\}} |v_j|^2 \, d\mu_w
\ge
\int_{B_j} |v_j|^2 \, d\mu_w
=
\|v_j\|_F^2
\ge
\lambda.
\end{equation}

Hence
\begin{equation}
\sup_{u \in \mathcal{B}_E} \int_{\{|z| > R\}} |u|^2 \, d\mu_w
\ge
\lambda.
\end{equation}

This holds for all $R > 0$, which contradicts the global tightness condition. Therefore $\mu_w(\mathbb{H}^{N+1}) < \infty$.
\end{proof}

\section{Main Embedding Theorems}\label{sec:Main Embedding Theorems}

The culmination of this study is presented in Section~\ref{sec:Main Embedding Theorems}, where the fundamental embedding theorems are formally stated and proved. Under the natural hypotheses of Strict Local Regularity and Uniform Local Comparability , Theorem~\ref{thm:Fundamental Characterization of Global Compact Embeddings} establishes a strict equivalence: the Rellich-Kondrachov embedding is compact if and only if the Sobolev unit ball is globally tight in the weighted Lebesgue space. This equivalence elegantly reduces the analytic difficulty of the embedding problem to verifiable geometric and measure-theoretic conditions, namely that the total mass is finite and the embedding admits uniform asymptotic mass-vanishing at infinity and at the singular boundary.

\begin{definition}[Strict Local Regularity]\label{def:Strict Local Regularity}
Let $\mathbb{H}^{N+1}$ be equipped with the subspace topology inherited from $\mathbb{R}^{N+1}$. Let $\mu_w$ be a Radon measure defined by $d \mu_w=w(z) d z$. The measure $\mu_w$ satisfies the Strict Local Regularity if, for every relatively compact, Lipschitz-regular open subset $\Omega \subset \mathbb{H}^{N+1}$ that is strictly separated from the boundary $\partial \mathbb{H}^{N+1}$ (i.e., $\inf _{(x, y) \in \Omega} y>0$ ), there exist constants $\lambda_{\Omega}, \Lambda_{\Omega} \in \mathbb{R}$ such that $0<\lambda_{\Omega} \leq \Lambda_{\Omega}<\infty$, satisfying:
$\lambda_{\Omega} \leq w(z) \leq \Lambda_{\Omega} \quad$ for Lebesgue-almost every $z \in \Omega$.
\end{definition}

\begin{proposition}[Isomorphism of Weighted and Unweighted Sobolev Spaces]\label{prop:Isomorphism of Weighted and Unweighted Sobolev Spaces}
Let $\mu_w$ satisfy Strict Local Regularity. Let $\Omega \subset \mathbb{H}^{N+1}$ be a relatively compact, Lipschitz-regular open subset that is strictly separated from the boundary $\partial \mathbb{H}^{N+1}$. Then the identity mapping
\begin{equation}
\operatorname{Id} : H^1_{\mu_w}(\Omega) \to H^1(\Omega, dz)
\end{equation}
is an isomorphism of Hilbert spaces, and the induced norms are strictly equivalent.
\end{proposition}

\begin{proof}
Let $V$ be the vector space of weakly differentiable functions on $\Omega$. Equip $V$ with two Hilbert space structures given by the norms $\|\cdot\|_{H^1(\Omega, dz)}$ and $\|\cdot\|_{H^1_{\mu_w}(\Omega)}$. Denote these spaces by
\begin{equation}
E_1 = (V, \|\cdot\|_{H^1(\Omega, dz)}), \qquad E_2 = (V, \|\cdot\|_{H^1_{\mu_w}(\Omega)}).
\end{equation}

Let $u \in V$. The unweighted Sobolev norm is given by
\begin{equation}
\|u\|_{H^1(\Omega, dz)}^2 = \int_\Omega \bigl(|u(z)|^2 + |\nabla u(z)|^2\bigr) \, dz,
\end{equation}
and the weighted Sobolev norm is given by
\begin{equation}
\|u\|_{H^1_{\mu_w}(\Omega)}^2 = \int_\Omega \bigl(|u(z)|^2 + |\nabla u(z)|^2\bigr) w(z) \, dz.
\end{equation}

By Strict Local Regularity, since $\Omega$ is relatively compact and strictly separated from $\partial \mathbb{H}^{N+1}$, there exist constants $\lambda_\Omega, \Lambda_\Omega > 0$ such that for Lebesgue-almost every $z \in \Omega$,
\begin{equation}
\lambda_\Omega \le w(z) \le \Lambda_\Omega.
\end{equation}
By monotonicity of the integral,
\begin{equation}
\int_\Omega \lambda_\Omega \bigl(|u|^2 + |\nabla u|^2\bigr) \, dz
\le
\int_\Omega \bigl(|u|^2 + |\nabla u|^2\bigr) w(z) \, dz
\le
\int_\Omega \Lambda_\Omega \bigl(|u|^2 + |\nabla u|^2\bigr) \, dz.
\end{equation}
Factoring out constants gives
\begin{equation}
\lambda_\Omega \|u\|_{H^1(\Omega, dz)}^2
\le
\|u\|_{H^1_{\mu_w}(\Omega)}^2
\le
\Lambda_\Omega \|u\|_{H^1(\Omega, dz)}^2.
\end{equation}
Since $\lambda_\Omega > 0$, taking square roots yields
\begin{equation}
\sqrt{\lambda_\Omega} \, \|u\|_{H^1(\Omega, dz)}
\le
\|u\|_{H^1_{\mu_w}(\Omega)}
\le
\sqrt{\Lambda_\Omega} \, \|u\|_{H^1(\Omega, dz)}.
\end{equation}

Define $\operatorname{Id} : E_2 \to E_1$ by $\operatorname{Id}(u) = u$. Then
\begin{equation}
\|\operatorname{Id}(u)\|_{E_1} \le \frac{1}{\sqrt{\lambda_\Omega}} \|u\|_{E_2},
\end{equation}
so $\operatorname{Id}$ is bounded. Similarly,
\begin{equation}
\|\operatorname{Id}^{-1}(u)\|_{E_2} \le \sqrt{\Lambda_\Omega} \|u\|_{E_1},
\end{equation}
so the inverse is bounded.

Since $E_1$ and $E_2$ share the same underlying vector space $V$, the mapping is bijective. A bijective bounded linear operator with bounded inverse between Hilbert spaces is an isomorphism. Hence $\operatorname{Id}$ is an isomorphism of Hilbert spaces.
\end{proof}

\begin{definition}[Smooth Exhaustion Projection Operators]
Let $\mathbb{H}^{N+1}$ be the open half-space. Let $\left(\Omega_m\right)_{m \in \mathbb{N}}$ be an admissible smooth exhaustion of $\mathbb{H}^{N+1}$ satisfying the definition. For each $m \in \mathbb{N}$, define the linear operator $P_m$ :
$F \rightarrow F$ by:

$$
P_m(f)=\chi_{\Omega_m} f
$$

where $F=L_{\mu_w}^2\left(\mathbb{H}^{N+1}\right)$.
\end{definition}

\begin{lemma}[Orthogonal Projection Properties of smooth Exhaustion]\label{lem:Orthogonal Projection Properties of smooth Exhaustion}
Let $F = L_{\mu_w}^2(\mathbb{H}^{N+1})$ be a Hilbert space where $\mu_w$ is a Radon measure. Let $(\Omega_m)_{m \in \mathbb{N}}$ be an admissible smooth exhaustion of $\mathbb{H}^{N+1}$. For each $m \in \mathbb{N}$, define the linear operator $P_m : F \to F$ by
\begin{equation}
P_m(f) = \chi_{\Omega_m} f.
\end{equation}
Then the sequence $(P_m)_{m \in \mathbb{N}}$ satisfies the following properties:

\begin{enumerate}
    \item For all $m \in \mathbb{N}$, $P_m$ is a well-defined orthogonal projection on $F$, and
\begin{equation}
\sup_m \|P_m\|_{\mathcal{L}(F)} \le 1.
\end{equation}
    \item For every $f \in F$,
\begin{equation}
\lim_{m \to \infty} \|f - P_m(f)\|_F = 0.
\end{equation}
\end{enumerate}
\end{lemma}

\begin{proof}
Fix $m \in \mathbb{N}$. By definition, $\Omega_m$ is an open subset of $\mathbb{H}^{N+1}$. Since $\mu_w$ is a Radon measure on $\mathbb{H}^{N+1}$, every open set is Borel and hence $\mu_w$-measurable. Therefore the indicator function $\chi_{\Omega_m}$ is $\mu_w$-measurable.

Let $f \in F$. Then $P_m(f) = \chi_{\Omega_m} f$ is measurable, and
\begin{equation}
\|P_m(f)\|_F^2 = \int_{\mathbb{H}^{N+1}} |\chi_{\Omega_m} f|^2 \, d\mu_w = \int_{\Omega_m} |f|^2 \, d\mu_w.
\end{equation}
Since $\Omega_m \subset \mathbb{H}^{N+1}$, monotonicity of the integral yields
\begin{equation}
\int_{\Omega_m} |f|^2 \, d\mu_w \le \int_{\mathbb{H}^{N+1}} |f|^2 \, d\mu_w = \|f\|_F^2.
\end{equation}
Hence $\|P_m(f)\|_F \le \|f\|_F$, so $P_m(f) \in F$ and $\|P_m\|_{\mathcal{L}(F)} \le 1$. Moreover,
\begin{equation}
P_m(P_m(f)) = \chi_{\Omega_m}^2 f = \chi_{\Omega_m} f = P_m(f),
\end{equation}
so $P_m$ is idempotent and therefore a projection. Taking the supremum over $m$ gives
\begin{equation}
\sup_m \|P_m\|_{\mathcal{L}(F)} \le 1.
\end{equation}

Now fix $f \in F$ and consider
\begin{equation}
\|f - P_m(f)\|_F^2 = \int_{\mathbb{H}^{N+1}} |f - \chi_{\Omega_m} f|^2 \, d\mu_w = \int_{\mathbb{H}^{N+1}} |f|^2 \chi_{\mathbb{H}^{N+1} \setminus \Omega_m} \, d\mu_w.
\end{equation}
Define
\begin{equation}
g_m(z) = |f(z)|^2 \chi_{\mathbb{H}^{N+1} \setminus \Omega_m}(z).
\end{equation}
Since $(\Omega_m)$ is an exhaustion of $\mathbb{H}^{N+1}$, we have
\begin{equation}
\bigcup_{m=1}^\infty \Omega_m = \mathbb{H}^{N+1}.
\end{equation}
Thus for every $z \in \mathbb{H}^{N+1}$, there exists $M \in \mathbb{N}$ such that $z \in \Omega_m$ for all $m \ge M$, which implies
\begin{equation}
\chi_{\mathbb{H}^{N+1} \setminus \Omega_m}(z) = 0 \quad \text{for all } m \ge M.
\end{equation}
Hence $g_m(z) \to 0$ pointwise for all $z \in \mathbb{H}^{N+1}$.

Moreover, for all $m$ and all $z$,
\begin{equation}
g_m(z) \le |f(z)|^2.
\end{equation}
Define $G(z) = |f(z)|^2$. Since $f \in F$, we have
\begin{equation}
\int_{\mathbb{H}^{N+1}} G \, d\mu_w = \|f\|_F^2 < \infty,
\end{equation}
so $G$ is integrable.

The sequence $(g_m)$ consists of measurable functions, converges pointwise to $0$, and is dominated by the integrable function $G$. By the Dominated Convergence Theorem,
\begin{equation}
\lim_{m \to \infty} \int_{\mathbb{H}^{N+1}} g_m \, d\mu_w = 0.
\end{equation}
Therefore,
\begin{equation}
\lim_{m \to \infty} \|f - P_m(f)\|_F^2 = 0,
\end{equation}
which implies
\begin{equation}
\lim_{m \to \infty} \|f - P_m(f)\|_F = 0.
\end{equation}
\end{proof}

\begin{proposition}[Global Compactness from Local Compactness and Global Tightness]\label{prop:Global Compactness from Local Compactness and Global Tightness}
Let $E = H_{\mu_w}^1(\mathbb{H}^{N+1})$ and $F = L_{\mu_w}^2(\mathbb{H}^{N+1})$. Let $\iota : E \to F$ be the canonical linear embedding. Let $(\Omega_m)_{m \in \mathbb{N}}$ be an admissible smooth exhaustion with associated orthogonal projections
\begin{equation}
P_m(f) = \chi_{\Omega_m} f.
\end{equation}
If the unit ball $\mathcal{B}_E$ is globally tight in $F$, then the embedding $\iota$ is compact.
\end{proposition}

\begin{proof}
Fix $m \in \mathbb{N}$. By definition, $\Omega_m$ is bounded, strictly separated from $\{y = 0\}$, and has boundary of class $C^\infty$. Hence $\Omega_m$ is Lipschitz-regular. By the strict local regularity condition (Proposition~\ref{prop:Isomorphism of Weighted and Unweighted Sobolev Spaces}, the weighted measure $d\mu_w$ is bounded above and below by strictly positive constants on $\Omega_m$. It follows that the weighted Sobolev and Lebesgue norms on $\Omega_m$ are equivalent to the unweighted norms. By the Rellich--Kondrachov theorem, the local embedding is compact, and therefore the operator $P_m \circ \iota : E \to F$ is compact for every $m \in \mathbb{N}$.

For each $k \in \mathbb{N}_{\ge 1}$, define
\begin{equation}
W_k = \{ z \in \mathbb{H}^{N+1} : |z| \le k \text{ and } y \ge 1/k \}.
\end{equation}
Since $W_k$ is closed and bounded in $\mathbb{R}^{N+1}$, it is compact. The family $(\Omega_m)_{m \in \mathbb{N}}$ consists of open sets whose union is $\mathbb{H}^{N+1}$, so it forms an open cover of $W_k$. By compactness, there exists a finite subcover, which is contained in a single $\Omega_{m(k)}$ due to the nested property of the exhaustion. Hence for each $k$ there exists $m(k) \in \mathbb{N}$ such that
\begin{equation}
W_k \subset \Omega_{m(k)}.
\end{equation}
Since $W_k \subset W_{k+1}$, the map $k \mapsto m(k)$ may be chosen increasing.

Let $m \ge m(k)$. Then $\Omega_m \supset \Omega_{m(k)} \supset W_k$, and taking complements gives
\begin{equation}
\Omega_m^c \subset W_k^c.
\end{equation}
By definition of $W_k$,
\begin{equation}
W_k^c = \{ z \in \mathbb{H}^{N+1} : |z| > k \} \cup \{ z \in \mathbb{H}^{N+1} : 0 < y < 1/k \}.
\end{equation}

Let $u \in \mathcal{B}_E$. Then
\begin{equation}
\|(I - P_m)\iota(u)\|_F^2 = \int_{\Omega_m^c} |u|^2 \, d\mu_w.
\end{equation}
By monotonicity and subadditivity of the integral,
\begin{equation}
\int_{\Omega_m^c} |u|^2 \, d\mu_w \le \int_{W_k^c} |u|^2 \, d\mu_w
\end{equation}
and
\begin{equation}
\int_{W_k^c} |u|^2 \, d\mu_w \le \int_{\{|z| > k\}} |u|^2 \, d\mu_w + \int_{\{0 < y < 1/k\}} |u|^2 \, d\mu_w.
\end{equation}

Taking the supremum over $u \in \mathcal{B}_E$, we obtain
\begin{equation}
\sup_{u \in \mathcal{B}_E} \|(I - P_m)\iota(u)\|_F^2 \le \sup_{u \in \mathcal{B}_E} \int_{\{|z| > k\}} |u|^2 \, d\mu_w + \sup_{u \in \mathcal{B}_E} \int_{\{0 < y < 1/k\}} |u|^2 \, d\mu_w.
\end{equation}

By global tightness, the first term converges to $0$ as $k \to \infty$, and the second term also converges to $0$ as $k \to \infty$. Hence
\begin{equation}
\lim_{m \to \infty} \sup_{u \in \mathcal{B}_E} \|(I - P_m)\iota(u)\|_F = 0.
\end{equation}

Thus $P_m \circ \iota$ is compact for all $m$, and the remainder converges to zero uniformly on the unit ball. It follows that $\iota$ is compact.
\end{proof}

\begin{theorem}[Fundamental Characterization of Global Compact Embeddings]\label{thm:Fundamental Characterization of Global Compact Embeddings}
Let $\mathbb{H}^{N+1}$ be the open half-space. Let $\mu_w$ be a Radon measure defined by $d\mu_w = w(z)\,dz$. Let $E = H_{\mu_w}^1(\mathbb{H}^{N+1})$ and $F = L_{\mu_w}^2(\mathbb{H}^{N+1})$. Let $\iota : E \hookrightarrow F$ be the canonical continuous linear embedding, and let $\mathcal{B}_E$ denote the closed unit ball of $E$. Assume the measure space satisfies Strict Local Regularity (Definition~\ref{def:Strict Local Regularity}) and Uniform Local Comparability (Definition~\ref{def:uniform local comparability}. Then the following are equivalent:

\begin{enumerate}
\item The canonical embedding $\iota : E \hookrightarrow F$ is compact.

\item The unit ball $\mathcal{B}_E$ is globally tight in $F$.

\item The measure satisfies $\mu_w(\mathbb{H}^{N+1}) < \infty$, and the embedding admits uniform asymptotic mass-vanishing at infinity and at the singular boundary.
\end{enumerate}

\end{theorem}

\begin{proof}
Assume that $\mathcal{B}_E$ is globally tight in $F$. Let $(\Omega_m)_{m \in \mathbb{N}}$ be an admissible smooth exhaustion of $\mathbb{H}^{N+1}$. For each $m \in \mathbb{N}$, define the linear orthogonal projection $P_m : F \to F$ by
\begin{equation}
P_m(f) = \chi_{\Omega_m} f.
\end{equation}
By the Lemma~\ref{lem:Topological Exhaustion}, the sequence $(P_m)$ is uniformly bounded by $1$, and for every $f \in F$ one has
\begin{equation}
\|f - P_m(f)\|_F \to 0 \quad \text{as } m \to \infty.
\end{equation}

Fix $m \in \mathbb{N}$. By the definition of admissible smooth exhaustion, $\Omega_m$ is bounded, separated from $\{y = 0\}$, and has boundary of class $C^\infty$. By the lemma, $\Omega_m$ is Lipschitz-regular. The Strict Local Regularity condition implies that the weighted norms on $\Omega_m$ are equivalent to the unweighted norms. Therefore, by the classical Rellich--Kondrachov theorem, the localized embedding
\begin{equation}
P_m \circ \iota : E \to F
\end{equation}
is compact for every $m \in \mathbb{N}$.

By the Proposition~\ref{prop:Global Compactness from Local Compactness and Global Tightness}, the orthogonal complement over the unit ball is bounded by the tail and boundary supremum integrals, and the global tightness hypothesis yields
\begin{equation}
\sup_{u \in \mathcal{B}_E} \|(I - P_m)\iota(u)\|_F \to 0
\quad \text{as } m \to \infty.
\end{equation}
Thus the hypotheses of the lemma are satisfied. It follows that the embedding $\iota$ is compact. This proves that (ii) implies (i).

Assume now that $\iota$ is compact. Then the image $\iota(\mathcal{B}_E)$ is relatively compact and therefore totally bounded in $F$. Let $\varepsilon > 0$ be given. There exists a finite $(\varepsilon/2)$-net
\begin{equation}
\mathcal{N} = \{f_1, \dots, f_k\} \subset F
\end{equation}
such that for every $u \in \mathcal{B}_E$ there exists some $f_j \in \mathcal{N}$ with
\begin{equation}
\|\iota(u) - f_j\|_F < \varepsilon/2.
\end{equation}

Fix $i \in \{1, \dots, k\}$. Since $f_i \in L_{\mu_w}^2(\mathbb{H}^{N+1})$, the function $|f_i|^2$ is $\mu_w$-integrable. Consider the measurable tail sets
\begin{equation}
E_R = \{ z \in \mathbb{H}^{N+1} : |z| > R \}.
\end{equation}
As $R \to \infty$, the sets $E_R$ decrease to the empty set. By continuity of the integral from above for integrable functions,
\begin{equation}
\lim_{R \to \infty} \int_{E_R} |f_i|^2 \, d\mu_w = 0.
\end{equation}
Since the set $\mathcal{N}$ is finite, there exists $R_0 > 0$ such that for every $R \ge R_0$ and every $i \in \{1, \dots, k\}$,
\begin{equation}
\|f_i\|_{F(E_R)} < \varepsilon/2.
\end{equation}

Let $u \in \mathcal{B}_E$, and choose $f_j \in \mathcal{N}$ such that $\|\iota(u) - f_j\|_F < \varepsilon/2$. Then
\begin{equation}
\|\iota(u)\|_{F(E_R)} \le \|\iota(u) - f_j\|_{F(E_R)} + \|f_j\|_{F(E_R)}.
\end{equation}
Since restriction decreases the norm,
\begin{equation}
\|\iota(u) - f_j\|_{F(E_R)} \le \|\iota(u) - f_j\|_F < \varepsilon/2.
\end{equation}
Hence for all $R \ge R_0$,
\begin{equation}
\|\iota(u)\|_{F(E_R)} < \varepsilon.
\end{equation}
Taking the supremum over $u \in \mathcal{B}_E$, we obtain
\begin{equation}
\sup_{u \in \mathcal{B}_E} \int_{\{|z| > R\}} |u|^2 \, d\mu_w \to 0
\quad \text{as } R \to \infty.
\end{equation}
This proves tail tightness.

The same argument applies to the boundary regions
\begin{equation}
B_\delta = \{ z \in \mathbb{H}^{N+1} : 0 < y < \delta \}.
\end{equation}
For each $f_i \in \mathcal{N}$, the sets $B_\delta$ decrease to the empty set as $\delta \to 0$, and therefore
\begin{equation}
\lim_{\delta \to 0} \int_{B_\delta} |f_i|^2 \, d\mu_w = 0.
\end{equation}
Since $\mathcal{N}$ is finite, there exists $\delta_0 > 0$ such that for every $0 < \delta \le \delta_0$ and every $i \in \{1, \dots, k\}$,
\begin{equation}
\|f_i\|_{F(B_\delta)} < \varepsilon/2.
\end{equation}
If $u \in \mathcal{B}_E$ and $f_j \in \mathcal{N}$ is chosen so that $\|\iota(u) - f_j\|_F < \varepsilon/2$, then
\begin{equation}
\|\iota(u)\|_{F(B_\delta)} \le \|\iota(u) - f_j\|_{F(B_\delta)} + \|f_j\|_{F(B_\delta)} < \varepsilon.
\end{equation}
Hence
\begin{equation}
\sup_{u \in \mathcal{B}_E} \int_{\{0 < y < \delta\}} |u|^2 \, d\mu_w \to 0
\quad \text{as } \delta \to 0.
\end{equation}
Thus $\mathcal{B}_E$ is globally tight in $F$. This proves that (i) implies (ii).

Assume next that (ii) holds. By the Lemma~\ref{prop:Global Compactness from Local Compactness and Global Tightness}, the global tightness of $\mathcal{B}_E$, together with the Uniform Local Comparability condition, implies
\begin{equation}
\mu_w(\mathbb{H}^{N+1}) < \infty.
\end{equation}
Moreover, the uniform asymptotic mass-vanishing at infinity and at the boundary are exactly the defining components of global tightness. Therefore (iii) holds.

Conversely, assume that (iii) holds. The asserted mass-vanishing conditions at infinity and at the singular boundary are precisely the requirements for tail tightness and boundary tightness. Hence $\mathcal{B}_E$ is globally tight in $F$, so (ii) holds.

Therefore (i), (ii), and (iii) are equivalent.
\end{proof}

\section{Conclusion}

In summary, we have established a definitive topological characterization of the Rellich-Kondrachov compact embedding on the half-space equipped with general singular weights $w(z) = y^c \phi(|z|)$. The geometric and analytic difficulties induced by both the unbounded nature of the domain and the boundary singularity are completely resolved by demonstrating the strict equivalence between global compactness and the uniform global tightness of the unit ball. The introduction of structural tools—namely, smooth exhaustion projections, localized fractional Hardy inequalities, and Lyapunov tail coercivity—yields a unified functional framework. This framework subsumes specific existing criteria for polynomial and Gaussian weights, providing the necessary foundation for future applications to degenerate elliptic equations and fractional extension problems.

\section{Conflict of Interest Statement}
The authors declare that they have no competing interests and receive no external funding.

\section{Data Availability Statement}
Data sharing is not applicable to this article as no datasets were generated or analyzed during the current study.

\end{document}